\documentclass[11pt]{article}
\usepackage[dvips, hmargin=1.5in, vmargin={1.5in, 2in}]{geometry}
\usepackage{amssymb, amsmath, amsthm, enumerate, stmaryrd}
\usepackage[latin1]{inputenc}
\usepackage[dvips, arrow, matrix, tips, curve]{xy}
\SelectTips{cm}{10}

 \DeclareMathOperator{\colim}{colim}

\newcommand{\cat}[1]{\mathbf{#1}}
\newcommand{\op}{\mathrm{op}}
\newcommand{\id}{\mathrm{id}}
\newcommand{\thg}{{\mathord{\text{--}}}}

\newcommand{\set}[2]{\left\{\,#1 \mid #2\,\right\}}
\newcommand{\spn}[1]{{\left<{#1}\right>}}
\newcommand{\elt}[1]{\left\llcorner{#1}\right\lrcorner}

\newcommand{\cd}[2][]{\vcenter{\hbox{\xymatrix#1{#2}}}}
\newcommand{\cdl}[2][]{\xymatrix@1#1{#2}}

\newcommand{\A}{{\mathcal A}}
\newcommand{\B}{{\mathcal B}}
\newcommand{\C}{{\mathcal C}}
\newcommand{\D}{{\mathcal D}}
\newcommand{\E}{{\mathcal E}}

\newcommand{\G}{{\mathcal G}}
\renewcommand{\H}{{\mathcal H}}
\newcommand{\I}{{\mathcal I}}
\newcommand{\J}{{\mathcal J}}
\newcommand{\K}{{\mathcal K}}
\renewcommand{\L}{{\mathcal L}}
\newcommand{\ELL}{{\mathcal L}}
\newcommand{\M}{{\mathcal M}}

\renewcommand{\P}{{\mathcal P}}

\newcommand{\R}{{\mathcal R}}
\newcommand{\T}{{\mathcal T}}

\newcommand{\V}{{\mathcal V}}

\newcommand{\xtor}[1]{\cdl[@1]{{} \ar[r]|-{\object@{|}}^{#1} & {}}}

\makeatletter

\def\hookleftarrowfill@{\arrowfill@\leftarrow\relbar{\relbar\joinrel\rhook}}
\def\twoheadleftarrowfill@{\arrowfill@\twoheadleftarrow\relbar\relbar}
\def\leftbararrowfill@{\arrowdoublefill@{\leftarrow\mkern-5mu}\relbar\mapstochar\relbar\relbar}
\def\Leftbararrowfill@{\arrowdoublefill@{\Leftarrow\mkern-2mu}\Relbar\Mapstochar\Relbar\Relbar}
\def\leftringarrowfill@{\arrowdoublefill@{\leftarrow\mkern-3mu}\relbar{\mkern-3mu\circ\mkern-2mu}\relbar\relbar}
\def\lefttriarrowfill@{\arrowfill@{\mathrel\triangleleft\mkern0.5mu\joinrel\relbar}\relbar\relbar}
\def\Lefttriarrowfill@{\arrowfill@{\mathrel\triangleleft\mkern1mu\joinrel\Relbar}\Relbar\Relbar}

\def\hookrightarrowfill@{\arrowfill@{\lhook\joinrel\relbar}\relbar\rightarrow}
\def\twoheadrightarrowfill@{\arrowfill@\relbar\relbar\twoheadrightarrow}
\def\rightbararrowfill@{\arrowdoublefill@{\relbar\mkern-0.5mu}\relbar\mapstochar\relbar\rightarrow}
\def\Rightbararrowfill@{\arrowdoublefill@{\Relbar\mkern-2mu}\Relbar\Mapstochar\Relbar\Rightarrow}
\def\rightringarrowfill@{\arrowdoublefill@\relbar\relbar{\mkern-2mu\circ\mkern-3mu}\relbar{\mkern-3mu\rightarrow}}
\def\righttriarrowfill@{\arrowfill@\relbar\relbar{\relbar\joinrel\mkern0.5mu\mathrel\triangleright}}
\def\Righttriarrowfill@{\arrowfill@\Relbar\Relbar{\Relbar\joinrel\mkern1mu\mathrel\triangleright}}

\def\leftrightarrowfill@{\arrowfill@\leftarrow\relbar\rightarrow}
\def\mapstofill@{\arrowfill@{\mapstochar\relbar}\relbar\rightarrow}

\renewcommand*\xleftarrow[2][]{\ext@arrow 20{20}0\leftarrowfill@{#1}{#2}}
\providecommand*\xLeftarrow[2][]{\ext@arrow 60{22}0{\Leftarrowfill@}{#1}{#2}}
\providecommand*\xhookleftarrow[2][]{\ext@arrow 10{20}0\hookleftarrowfill@{#1}{#2}}
\providecommand*\xtwoheadleftarrow[2][]{\ext@arrow 60{20}0\twoheadleftarrowfill@{#1}{#2}}
\providecommand*\xleftbararrow[2][]{\ext@arrow 10{22}0\leftbararrowfill@{#1}{#2}}
\providecommand*\xLeftbararrow[2][]{\ext@arrow 50{24}0\Leftbararrowfill@{#1}{#2}}
\providecommand*\xleftringarrow[2][]{\ext@arrow 10{26}0\leftringarrowfill@{#1}{#2}}
\providecommand*\xlefttriarrow[2][]{\ext@arrow 80{24}0\lefttriarrowfill@{#1}{#2}}
\providecommand*\xLefttriarrow[2][]{\ext@arrow 80{24}0\Lefttriarrowfill@{#1}{#2}}

\renewcommand*\xrightarrow[2][]{\ext@arrow 01{20}0\rightarrowfill@{#1}{#2}}
\providecommand*\xRightarrow[2][]{\ext@arrow 04{22}0{\Rightarrowfill@}{#1}{#2}}
\providecommand*\xhookrightarrow[2][]{\ext@arrow 00{20}0\hookrightarrowfill@{#1}{#2}}
\providecommand*\xtwoheadrightarrow[2][]{\ext@arrow 03{20}0\twoheadrightarrowfill@{#1}{#2}}
\providecommand*\xrightbararrow[2][]{\ext@arrow 01{22}0\rightbararrowfill@{#1}{#2}}
\providecommand*\xRightbararrow[2][]{\ext@arrow 04{24}0\Rightbararrowfill@{#1}{#2}}
\providecommand*\xrightringarrow[2][]{\ext@arrow 01{26}0\rightringarrowfill@{#1}{#2}}
\providecommand*\xrighttriarrow[2][]{\ext@arrow 07{24}0\righttriarrowfill@{#1}{#2}}
\providecommand*\xRighttriarrow[2][]{\ext@arrow 07{24}0\Righttriarrowfill@{#1}{#2}}

\providecommand*\xmapsto[2][]{\ext@arrow 01{20}0\mapstofill@{#1}{#2}}
\providecommand*\xleftrightarrow[2][]{\ext@arrow 10{22}0\leftrightarrowfill@{#1}{#2}}
\providecommand*\xLeftrightarrow[2][]{\ext@arrow 10{27}0{\Leftrightarrowfill@}{#1}{#2}}

\makeatother

\newcommand{\twocong}[2][0.5]{\ar@{}[#2] \save ?(#1)*{\cong}\restore}
\newcommand{\twoeq}[2][0.5]{\ar@{}[#2] \save ?(#1)*{=}\restore}
\newcommand{\rtwocell}[3][0.5]{\ar@{}[#2] \ar@{=>}?(#1)+/l 0.2cm/;?(#1)+/r 0.2cm/^{#3}}
\newcommand{\ltwocell}[3][0.5]{\ar@{}[#2] \ar@{=>}?(#1)+/r 0.2cm/;?(#1)+/l 0.2cm/^{#3}}
\newcommand{\ltwocello}[3][0.5]{\ar@{}[#2] \ar@{=>}?(#1)+/r 0.2cm/;?(#1)+/l 0.2cm/_{#3}}
\newcommand{\dtwocell}[3][0.5]{\ar@{}[#2] \ar@{=>}?(#1)+/u  0.2cm/;?(#1)+/d 0.2cm/^{#3}}
\newcommand{\dthreecell}[3][0.5]{\ar@{}[#2] \ar@3{->}?(#1)+/u  0.2cm/;?(#1)+/d 0.2cm/^{#3}}
\newcommand{\utwocell}[3][0.5]{\ar@{}[#2] \ar@{=>}?(#1)+/d 0.2cm/;?(#1)+/u 0.2cm/_{#3}}
\newcommand{\dtwocelltarg}[3][0.5]{\ar@{}#2 \ar@{=>}?(#1)+/u  0.2cm/;?(#1)+/d 0.2cm/^{#3}}
\newcommand{\utwocelltarg}[3][0.5]{\ar@{}#2 \ar@{=>}?(#1)+/d  0.2cm/;?(#1)+/u 0.2cm/_{#3}}

\newcommand{\pushoutcorner}[1][dr]{\save*!/#1-1.2pc/#1:(-1,1)@^{|-}\restore}

\newdir{ >}{{}*!/-5pt/\dir{>}}

\makeatletter
\newcommand{\pgph}{\paragraph{\!}}
\renewcommand{\paragraph}{\@startsection
{paragraph}%
{3}%
{0mm}%
{-\baselineskip}%
{-0.4em plus 0.2em minus 0.2em}%
{\normalfont\normalsize\bfseries}}%
\makeatother

\numberwithin{equation}{section} \numberwithin{paragraph}{section}

\newenvironment{Thm}{\paragraph{Theorem:}\em}{\vskip\baselineskip}
\newenvironment{Defn}{\paragraph{Definition:}}{\vskip\baselineskip}

\newenvironment{Prop}{\paragraph{Proposition:}\em}{\vskip\baselineskip}
\newenvironment{Propstar}{\paragraph{Proposition$^\ast$:}\em}{\vskip\baselineskip}

\newenvironment{Rk}{\paragraph{Remark:}}{\vskip\baselineskip}
\newenvironment{Exs}{\paragraph{Examples:}}{\vskip\baselineskip}
\newenvironment{Ex}{\paragraph{Example:}}{\vskip\baselineskip}

\makeatletter \@namedef{itemize*}{\itemize\parsep\z@ \parskip\z@}
\@namedef{enditemize*}{\enditemize} \@namedef{enumerate*}{\enumerate\parsep\z@
\parskip\z@} \@namedef{endenumerate*}{\endenumerate}
\def\Pr@@f{\subsubsection*{\textbf{Proof}}}
\def\pr@@f[#1]{\subsubsection*{{\textbf{Proof}} #1}}
\makeatother

\DeclareMathOperator{\Lan}{Lan}

\makeatletter \@namedef{itemize*}{\itemize\parsep\z@ \parskip\z@}
\@namedef{enditemize*}{\enditemize} \@namedef{enumerate*}{\enumerate\parsep\z@
\parskip\z@} \@namedef{endenumerate*}{\endenumerate} \makeatother

\begin{document}

\title{Understanding the small object argument}
\author{Richard Garner\thanks{Supported by a Research Fellowship of St John's College, Cambridge and a Marie Curie Intra-European Fellowship, Project No.\ 040802.}\\Department of Mathematics, Uppsala University,\\Box 480, S-751 06 Uppsala, Sweden}
\maketitle
\begin{abstract}
The small object argument is a transfinite construction which, starting from a set of maps in
a category, generates a weak factorisation system on that category. As useful as it is, the
small object argument has some problematic aspects: it possesses no universal property; it
does not converge; and it does not seem to be related to other transfinite constructions
occurring in categorical algebra. In this paper, we give an ``algebraic'' refinement of the
small object argument, cast in terms of Grandis and Tholen's natural weak factorisation
systems, which rectifies each of these three deficiencies.
\end{abstract}

%GATHER{biblio.bib}
\newcommand{\inj}{\textsf{in}}
\newcommand{\catc}{\K}%\cat{Cat} / \Ar \C}
\newcommand{\comp}{\cat{Mon}(\K)}
\newcommand{\orth}{\mathop{\boxempty}}
\newcommand{\morth}{\mathop{\boxtimes}}
\newcommand{\lp}{\cat{Sq}}
\newcommand{\zig}{\cat{ZZ}}
\newcommand{\sorth}{\mathop{\bot}}
\newcommand{\lo}[1]{{}^{\orth} {#1}}
\newcommand{\ro}[1]{{#1}^{\orth}}
\newcommand{\mlo}[1]{{\vphantom{#1}}^{\morth} {#1}}
\newcommand{\mro}[1]{{#1}^{\morth}}
\newcommand{\Ar}[1]{{{#1}^\mathbf{2}}}
\newcommand{\wfs}{w.f.s.}
\newcommand{\nwfs}{n.w.f.s.}
\newcommand{\dom}{\mathrm{dom}}
\newcommand{\cod}{\mathrm{cod}}
\newcommand{\map}{\mathbf}
\newcommand{\Ll}{\mathsf L}
\newcommand{\Rr}{\mathsf R}
\newcommand{\Lmap}{\mathsf L\text-\cat{Map}}
\newcommand{\Rmap}{\mathsf R\text-\cat{Map}}
\newcommand{\Ff}{\cat{Ff}}
\let\L\undefined

\section{Introduction}
The concept of \emph{factorisation system} provides us with a way of viewing a
category $\C$ as a compositional product of two subcategories $\ELL$ and $\R$.
The two key ingredients are an axiom of \emph{factorisation}, which affirms
that any map of $\C$ may be written as a map of $\ELL$ followed by a map of
$\R$, and an axiom of \emph{orthogonality}, which assures us that this
decomposition is unique up to unique isomorphism. From these two basic axioms a
very rich theory can be developed, and a very useful one, since most categories
arising in mathematical practice will admit at least a few different
factorisation systems.

However, in those mathematical areas where the primary objects of study are
themselves higher-dimensional entities~--~most notably, topology and higher
dimensional category theory~--~the notion of factorisation system is frequently
too strong, since we would like factorisations be unique, not up to
isomorphism, but up to something weaker. Thus in a 2-category, we might want
uniqueness up-to-equivalence; or in a category of topological spaces,
uniqueness up-to-homotopy.

The usual way of achieving this is to pass from factorisation systems to
\emph{weak factorisation systems}. The modifier ``weak'' has the familiar
effect of turning an assertion of unique existence into an assertion of mere
existence, here in respect to the diagonal fill-ins which are guaranteed to us
by the axiom of orthogonality.

From this, we would not necessarily expect the factorisations in a weak
factorisation system (henceforth w.f.s.) to be unique up to anything at all:
but remarkably, each weak factorisation system generates its own notion of
``equivalence'' which respect to which its factorisations \emph{are} unique.
The framework within which this is most readily expressed is that of Quillen's
\emph{model categories} \cite{Quillen:homotopical}, which consist in a clever
interaction of two w.f.s.'s on a category: but we can make do with a single
w.f.s., and for the purposes of this paper, we will.

Whilst in many respects, the theory of w.f.s.'s is similar to the theory of
factorisation systems (which we will henceforth call \emph{strong}
factorisation systems to avoid ambiguity), there are some puzzling aspects to
it: and notable amongst these is the manner in which one typically constructs a
w.f.s.

In the case of strong factorisation systems, there is a very elegant theory
which, given a sufficiently well-behaved category $\C$, can generate a strong
factorisation system from any set of maps $J \subset \C^\mathbf 2$. The
$\R$-maps will be the maps which are \emph{right orthogonal} to each of the
maps in $J$ (in a sense which we recall more precisely in Section
\ref{Sec:nwfs}); and the $\ELL$-maps, those which are left orthogonal to each
of the maps in $\R$. The key difficulty is how we should build the
factorisations, and for this we are able bring to bear a well-established body
of knowledge concerning transfinite constructions in categories, on which the
definitive word is \cite{Ke80}.

There is a corresponding theory for weak factorisation systems. Again, we
suppose ourselves given a well-behaved $\C$ and a set of maps $J$, but this
time we take for $\R$ the class of maps \emph{weakly} right orthogonal to $J$,
and for $\ELL$, the class of maps \emph{weakly} left orthogonal to $\R$. To
obtain a weak factorisation system, we must also have factorisation of maps:
and for this, we apply a construction known as the \emph{small object
argument}, introduced by Quillen \cite{Quillen:homotopical}, and first given in
its full generality by Bousfield \cite{Bous}.

The problem lies in divining the precise nature of the small object argument.
It is certainly some kind of transfinite construction: but it is a transfinite
construction which does not converge, has no universal property, and does not
seem to be an instance of any other known transfinite construction.

In this paper, we present a modification of the small object argument which
rectifies each of these deficiencies: it is guaranteed to converge; the
factorisations it provides are freely generated by the set $J$, in a suitable
sense; and it may be construed as an instance of a familiar free monoid
construction.

To make this possible, we must adopt a rather different perspective on weak
factorisation systems. The definition of a w.f.s.\ specifies classes of maps
$\ELL$ and $\R$ together with axioms which affirm properties: that \emph{there
exist} factorisations, or that \emph{there exist} certain diagonal fill-ins.
But a key tenet of category theory is that anything we specify in terms of
properties should have an equally valid expression in terms of structure: and
in the case of w.f.s.'s, a suitable ``algebraic'' reformulation is given by
Tholen and Grandis' notion of \emph{natural weak factorisation system}
\cite{nwfs}.

The extra algebraicity provided by natural w.f.s.'s allows us a clearer view of
what is actually going on in the small object argument. We now have a
\emph{functor} from the category of natural w.f.s.'s on $\C$ into $\cat{CAT}$
which sends each natural w.f.s.\ to its category of $\ELL$-maps; and we can
factor this functor through $\cat{CAT} / \C^\mathbf 2$. We may view the
resultant functor $\cat{NWFS}(\C) \to \cat{CAT} / \C^\mathbf 2$ as being the
``semantics'' side of a syntax/semantics adjunction: for which the syntax side
is precisely our refinement of the small object argument.

Although all our arguments will be cast in terms of natural w.f.s.'s, we will
see that there are ramifications for plain \wfs's as well, since our refined
version of the small object argument can equally well be applied there, giving
rise to factorisations which are less redundant than the original argument, and
in many cases can be easily calculated by hand.

\textbf{Acknowledgements}. My foremost thanks go to the organisers of CT '07
for providing such a pleasant and stimulating environment within which to
present this material. Further thanks go to Clemens Berger, Eugenia Cheng, Jeff
Egger, André Hirschowitz, Martin Hyland, Joachim Kock, Mike Shulman, Carlos
Simpson, Walter Tholen, and members of the Stockholm-Uppsala Logic Seminar for
useful discussions and comments.

\section{Notions of factorisation system}
\label{Sec:nwfs} In this section, we describe in detail the various sorts of
factorisation system mentioned in the Introduction.

\pgph Most familiar is the notion of \emph{strong factorisation system} $(\ELL,
\R)$ on a category $\C$, introduced by Freyd and Kelly in
\cite{FreydKelly:continuous}. This is given by two classes of maps $\ELL$ and
$\R$ in $\C$ which are each closed under composition with isomorphisms, and
which satisfy the axioms of
\begin{description}
\item[(factorisation)] Every map $e \colon X \to Y$ in $\C$ can be written
    as $e = gf$, where $f \in \ELL$ and $g \in \R$; and
\item[(orthogonality)] $f \mathbin\bot g$ for all $f \in \ELL$ and $g \in
    \R$, where $f \mathbin \bot g$ means that for every commutative square
\begin{equation}
\label{fillinsquare}
    \cd{
        A
            \ar[r]^-h
            \ar[d]_f &
        C
            \ar[d]^g \\
        B
            \ar[r]_-k &
        D
    }
\end{equation}
in $\C$, there is a unique map $j \colon B \to C$ such that $gj = k$ and
$jf = h$.
\end{description}
Instead of writing $f \mathbin \bot g$, we may also say that $f$ is \emph{left
orthogonal} to $g$ or that $g$ is \emph{right orthogonal} to $f$; moreover,
given a class $\A$ of maps in $\C$, we write
\begin{equation*}
    {}^\bot \A = \set{g \in \C^\mathbf 2}{f \mathbin \bot g \text{ for all $g \in \A$}}
\quad\text{and} \quad
    \A^\bot = \set{f \in \C^\mathbf 2}{f \mathbin \bot g \text{ for all $f \in \A$}}\text;
\end{equation*}
and this sets up a Galois connection on the collection of all classes of maps
in $\C$. In a strong factorisation system, we have $\R = \ELL^\bot$ and $\ELL =
{}^\bot \R$, so that the classes $\ELL$ and $\R$ determine each other.

\pgph \label{defwfs} We arrive at the notion of a \emph{weak factorisation
system} \cite{Bous} by making two alterations to the above definition. One is
minor: we require that $\ELL$ and $\R$ are closed under retracts in the arrow
category $\C^\mathbf 2$, rather than merely closed under isomorphism. The other
is more far-reaching: we replace the orthogonality condition with
\begin{description}
\item[(weak orthogonality)] $f \mathbin\pitchfork g$ for all $f \in \ELL$
    and $g \in \R$, where $f \mathbin\pitchfork g$ means that for every
    commutative square as in \eqref{fillinsquare}, there exists a (not
    necessarily unique) fill-in $j \colon B \to C$ such that $gj = k$ and
    $jf = h$.
\end{description}

\noindent We now have a Galois connection ${}^{\pitchfork}(\ ) \mathbin \dashv
(\ )^\pitchfork$; and again, the classes $\ELL$ and $\R$ of a w.f.s.\ determine
each other by the equations $\ELL = {}^{\pitchfork}\R$ and $\R =
\ELL^\pitchfork$. However, the classes $\ELL$ and $\R$ need not determine the
factorisation of a map, even up to isomorphism, as the following examples show:

\begin{Exs}\label{exs1}\hfill
\begin{enumerate}[(i)]
\item (Epi, Mono) is a strong factorisation system on $\cat{Set}$; but
    (Mono, Epi) is a weak factorisation system. For the latter, there are
    two natural choices of factorisation for a map $f \colon X \to Y$: the
    \emph{graph} factorisation which goes via $X \times Y$; and the
    \emph{cograph} factorisation, which goes through $X + Y$.
\item There is a weak factorisation system on $\cat{Cat}$ given by
    (injective equivalences, isofibrations). A \emph{injective equivalence}
    is a functor which is both injective on objects and an equivalence of
    categories; whilst an \emph{isofibration} is a functor along which all
    isomorphisms have liftings.
\item There is a weak factorisation system (anodyne extensions, Kan
    fibrations) on $\cat{SSet} = [\Delta^\op, \cat{Set}]$, the category of
    simplicial sets. The Kan fibrations are easy to describe: they are
    precisely the maps which are weakly right orthogonal to the set of
    \emph{horn inclusions} $\Lambda^k[n] \to \Delta[n]$. The anodyne
    extensions are the class of maps weakly left orthogonal to all Kan
    fibrations; more explicitly, they are obtained by closing the set of
    horn inclusions under countable composition, cobase change, coproduct
    and retract.
\end{enumerate}
\end{Exs}

\pgph As we mentioned in the Introduction, Grandis and Tholen's \emph{natural
weak factorisation systems} \cite{nwfs} provide an algebraisation of the notion
of weak factorisation system. In order to motivate the definition, we first
give a similar algebraisation of the notion of strong factorisation system.

So suppose that we are given a strong factorisation system $(\ELL, \R)$ on a
category $\C$, together with for each map of $\C$ a choice of factorisation:
\begin{equation*}
    X \xrightarrow f Y \quad \mapsto \quad X \xrightarrow{\lambda_f} Kf \xrightarrow{\rho_f} Y\text,
\end{equation*}
where $\lambda_f \in \ELL$ and $\rho_f \in \R$. It follows  from the
orthogonality property that this assignation may be extended in a unique way to
a \emph{functorial factorisation}: which is to say, a functor $F \colon
\C^\mathbf 2 \to \C^\mathbf 3$ (where $\mathbf 2$ and $\mathbf 3$ are the
ordinals $(0 \leqslant 1)$ and $(0 \leqslant 1 \leqslant 2)$ respectively)
which splits the ``face map'' $d_1 \colon \C^\mathbf 3 \to \C^\mathbf 2$ given
by
\begin{equation*}
    d_1(X \xrightarrow f Y \xrightarrow g Z) = (X \xrightarrow{gf} Z)\text.
\end{equation*}
\pgph This face map is induced by the functor $\delta_1 \colon \mathbf 2 \to
\mathbf 3$ picking out the unique arrow $0 \to 2$: we have $d_1 =
\C^{\delta_1}$. There are two other functors $\delta_0, \delta_2 \colon \mathbf
2 \to \mathbf 3$, and homming these into $\C$ induces two further face maps
$d_0, d_2 \colon \C^\mathbf 3 \to \C^\mathbf 2$, with
\begin{equation*}
    d_0(X \xrightarrow f Y \xrightarrow g Z) = (Y \xrightarrow{g} Z)\quad \text{and}\quad
    d_2(X \xrightarrow f Y \xrightarrow g Z) = (X \xrightarrow{f} Y)\text.
\end{equation*}
Postcomposing our functorial factorisation $F \colon \C^\mathbf 2 \to
\C^\mathbf 3$ with these induces functors $L, R \colon \C^\mathbf 2 \to
\C^\mathbf 2$, which send an object $f$ of $\C^\mathbf 2$ to $\lambda_f$ and
$\rho_f$ respectively.

\pgph \label{furtherstructure} There is further structure in $\cat{Cat}(\mathbf
2, \mathbf 3)$ which we can make use of: we have natural transformations
$\gamma_{2, 1} \colon \delta_2 \Rightarrow \delta_1$ and $\gamma_{1, 0} \colon
\delta_1 \Rightarrow \delta_0$, and by homming these into $\C$, we induce
natural transformations $c_{2, 1} \colon d_2 \Rightarrow d_1$ and $c_{1, 0}
\colon d_1 \Rightarrow d_0$. Postcomposing $F \colon \C^\mathbf 2 \to
\C^\mathbf 3$ with these now gives us natural transformations $\Phi \colon L
\Rightarrow \id_{\C^\mathbf 2}$ and $\Lambda \colon \id_{\C^\mathbf 2}
\Rightarrow R$ with components
\[
\Phi_f =
    \cd[@1]{
        X \ar[d]_{\lambda_f} \ar[r]^{\id_X} &
        X \ar[d]^{f} \\
        Kf \ar[r]_{\rho_f} &
        Y
    }
\qquad \text{and} \qquad \Lambda_f  =
    \cd[@1]{
        X \ar[d]_{f} \ar[r]^{\lambda_f} &
        Kf \ar[d]^{\rho_f} \\
        Y \ar[r]_{\id_Y} &
        Y\text.
    }
\]
\pgph Now, because $F \colon \C^\mathbf 2 \to \C^\mathbf 3$ arose from a strong
factorisation system, the corresponding $\Lambda : \id_{\C^\mathbf 2}
\Rightarrow R$ will provide the unit for a reflection of $\C^\mathbf 2$ into
the full subcategory of $\C^\mathbf 2$ spanned by the $\R$-maps. To see this,
consider a morphism
\begin{equation*}
    \cd[@1]{
        X \ar[d]_{f} \ar[r]^{h} &
        W \ar[d]^{g} \\
        Y \ar[r]_{k} &
        Z
    }
\end{equation*}
from $f$ to $g$ in $\C^\mathbf 2$, with $g$ an $\R$-map. Then applying
orthogonality to the square
\begin{equation*}
        \cd{
        X \ar[d]_{\lambda_f} \ar[r]^{h} &
        W \ar[d]^{g} \\
        Kf \ar[r]_{k . \rho_f} &
        Z
    }
\end{equation*}
we obtain a map $j \colon Kf \to W$ making both triangles commute; and so the
map $(h, k) \colon f \to g$ factors uniquely through $\Lambda_f$ as
\begin{equation*}
    f \xrightarrow{\Lambda_f} \rho_f \xrightarrow{(j, k)} g\text.
\end{equation*}

Thus the subcategory spanned by the $\R$-maps is a full, replete, reflective
subcategory of $\C^\mathbf 2$, via the reflector $\Lambda \colon
\id_{\C^\mathbf 2} \Rightarrow R$, and so $(R, \Lambda)$ extends uniquely to an
idempotent monad $\Rr = (R, \Lambda, \Pi)$ whose category of $\Rr$-algebras may
be identified with this subcategory. Dually, the pair $(L, \Phi)$ may be
extended uniquely to an idempotent comonad $\Ll = (L, \Phi, \Sigma)$ whose
category of coalgebras is isomorphic to the full subcategory of $\C^\mathbf 2$
spanned by the $\ELL$-maps. Thus we have proved:

\begin{Prop}\label{swfs}
There is a bijective correspondence between strong factorisation systems
$(\ELL, \R)$ on a category $\C$ for which a choice of factorisation for every
map has been made, and functorial factorisations $F \colon \C^\mathbf 2 \to
\C^\mathbf 3$ for which the corresponding pointed endofunctor $(R, \Lambda)$
underlies an idempotent monad and the corresponding copointed endofunctor $(L,
\Phi)$ underlies an idempotent comonad.
\end{Prop}

The notion of natural weak factorisation system now arises by generalising the
situation of this Proposition in a very obvious way: by dropping the
requirement of idempotency.

\begin{Defn}\cite{nwfs}
A \emph{natural weak factorisation system} on a category $\C$ is given by a
functorial factorisation $F \colon \C^\mathbf 2 \to \C^\mathbf 3$, together
with an extension of the corresponding pointed endofunctor $(R, \Lambda)$ to a
monad $\Rr = (R, \Lambda, \Pi)$; and an extension of the corresponding
copointed endofunctor $(L, \Phi)$ to a comonad $\Ll = (L, \Phi, \Sigma)$.
\end{Defn}
Observe that we can reconstruct $F$ from $\Ll$ and $\Rr$, and thus we may speak
simply of a natural weak factorisation system $(\Ll, \Rr)$.
\begin{Exs}\label{exs2}
\begin{enumerate}[(i)]
\item There is a natural w.f.s.\ on $\cat{Set}$ whose underlying functorial
    factorisation is the graph factorisation of Examples \ref{exs1}(i):
    \begin{equation*}
        X \xrightarrow f Y \quad \mapsto \quad X \xrightarrow{\spn{id, f}} X \times Y \xrightarrow{\pi_2} Y\text.
    \end{equation*}
    Dually, there is a natural w.f.s.\ on $\cat{Set}$ which factors $f$
    through $X+Y$. These examples generalise to any category with products
    or coproducts, as the case may be.
\item There is a natural w.f.s.\ on $\cat{Cat}$ whose underlying functorial
    factorisation is given by
    \begin{equation*}
        \C \xrightarrow F \D \quad \mapsto \quad \C \xrightarrow{\lambda_F} \D \downarrow F \xrightarrow {\rho_F} \D\text,
    \end{equation*}
    where $\D \downarrow F$ is the \emph{comma category} whose objects are
    triples $(c, d, f \colon d \to Fc)$; $\lambda_F$ is the functor sending
    $c$ in $\C$ to $(\id \colon Fc \to Fc)$ in $\D \downarrow F$; and
    $\rho_F$ is the functor sending $(f \colon d \to Fc)$ in $\D \downarrow
    F$ to $d$. There are variations on this theme: we can replace $\D
    \downarrow F$ with the dual comma category $F \downarrow \D$; or with
    the \emph{iso-comma category} $\D \downarrow_{\cong} F$, which is the
    full subcategory of $\D \downarrow F$ whose objects are the invertible
    arrows. These examples generalise to any 2-category with comma objects.
\item By Proposition \ref{swfs}, any strong factorisation system on $\C$
    gives  rise to a natural weak factorisation system on $\C$.
\end{enumerate}
\end{Exs}

\pgph It is not immediately clear that a natural w.f.s.\ deserves the name of
weak factorisation system. To show that this is so, we must exhibit suitable
analogues of the axioms of factorisation and weak orthogonality; for which we
must first identify what the $\ELL$-maps and $\R$-maps are. Now, for a strong
factorisation system, we can reconstruct the $\ELL$-\ and $\R$-maps from the
associated comonad $\Ll$ and monad $\Rr$ as their respective coalgebras and
algebras; and thus it is natural to define:
\begin{Defn}
Let $(\Ll, \Rr)$ be a natural w.f.s.\ on $\C$. We write $\Lmap$ for the
category of $\Ll$-coalgebras, and call its objects \emph{$\Ll$-maps}; and write
$\Rmap$ for the category of $\Rr$-algebras and call its objects
\emph{$\Rr$-maps}.
\end{Defn}
Note that being an $\Ll$-\ or $\Rr$-map is structure on, and not a property of,
a map of $\C$.
\begin{Exs}\label{exs3}
\begin{enumerate}[(i)]
\item For the natural w.f.s.\ on $\cat{Set}$ which factors $f \colon X \to
    Y$ through $X + Y$, an $\Rr$-map structure on $g \colon C \to D$ is a
    splitting for $g$: that is, a morphism $g^\ast \colon Y \to X$ with
    $gg^\ast = \id_Y$. An $\Ll$-map structure on $f \colon A \to B$ exists
    just when $f$ is a monomorphism, and in this case is uniquely
    determined: thus the comonad $\Ll$ is ``property-like'', though not
    idempotent.
\item For the natural w.f.s.\ on $\cat{Cat}$ which factors $F \colon \C \to
    \D$ through $\D \downarrow F$, an $\Rr$-map is a \emph{split
    fibration}: that is, a Grothendieck fibration with chosen liftings that
    compose up strictly. An $\Ll$-map is, roughly speaking, an inclusion of
    a reflective subcategory: more precisely, an $\Ll$-map structure on a
    functor $F \colon \C \to \D$ is given by specifying a functor $F^\ast
    \colon \D \to \C$ and a natural transformation $\eta \colon 1_\D
    \Rightarrow FF^\ast$ satisfying $F^\ast F = 1_\D$, $F^\ast \eta =
    \id_{F^\ast}$ and $\eta F = \id_F$. For the n.w.f.s.\ which factors
    through $F \downarrow \D$ instead, the $\Rr$-algebras are split
    opfibrations and the $\Ll$-coalgebras, inclusions of coreflective
    subcategories; whilst if we factor through $\D \downarrow_{\cong} F$,
    then $\Rr$-algebras are split isofibrations, and $\Ll$-coalgebras are
    retract equivalences.
\item If we view a strong factorisation system $(\ELL, \R)$ on $\C$ as a
    natural w.f.s., then the $\Ll$-maps and $\Rr$-maps reduce to
    $\ELL$-maps and $\R$-maps. In this particular case, being an $\Ll$-\ or
    $\Rr$-map returns to being a mere property; and this is because the
    comonad $\Ll$ and monad $\Rr$ are idempotent.
\end{enumerate}
\end{Exs}
Further details on these examples may be found in \cite{nwfs}.

 \pgph \label{factdefn} With this definition of $\Ll$-map and $\Rr$-map, it is now
clear that any natural w.f.s.\ $(\Ll, \Rr)$ admits an axiom of factorisation:
given a map $f \colon C \to D$, we obtain an $\Ll$-map structure on $\lambda_f
\colon C \to Kf$ by applying the cofree functor $\C^\mathbf 2 \to \Lmap$, and
an $\Rr$-map structure on $\rho_f \colon Kf \to D$ by applying the free functor
$\C^\mathbf 2 \to \Rmap$.

\pgph \label{worthdefn} More interestingly, we also have an axiom of weak
orthogonality. To see this, suppose that we are given a square like
\eqref{fillinsquare} together with an $\Ll$-coalgebra structure on $f$ and an
$\Rr$-algebra structure on $g$. Thus we have a coaction morphism $e \colon f
\to Lf$ and an action morphism $m \colon Rg \to g$, which the (co)algebra
axioms force to be of the following forms:
\begin{equation*}
    e = \cd{
        A \ar[d]_f \ar[r]^{\id_A} & A \ar[d]^{\lambda_f} \\
        B \ar[r]_{s} & Kf
    } \qquad \text{and} \qquad
    m = \cd{
        Kf \ar[d]_{\rho_g} \ar[r]^{p} & C \ar[d]^{g} \\
        D \ar[r]_{\id_D} & D\text.
    }
\end{equation*}
Furthermore, we may view the square \eqref{fillinsquare} as a map $(h, k)
\colon f \to g$ in $\C^\mathbf 2$; and so applying the functorial factorisation
of $(\Ll, \Rr)$ yields an arrow $K(h, k) \colon Kf \to Kg$ in $\C$. We now
obtain a diagonal fill-in for \eqref{fillinsquare} as the composite
\begin{equation}
    B \xrightarrow s Kf \xrightarrow{K(h, f)} Kg \xrightarrow p C\text.
    \label{fillineq}
\end{equation}
Note that this fill-in is canonically determined by the $\Ll$-map structure on
$f$ and the $\Rr$-map structure on $g$. Indeed, it is reasonable to view an
$\Ll$-map structure as encoding a coherent choice of lifting opposite every
$\Rr$-map, and vice versa.

\begin{Ex}
Let us see how we obtain diagonal fill-ins for the natural w.f.s.\ on
$\cat{Cat}$ which factors $F \colon \C \to \D$ through $\D \downarrow F$. We
suppose ourselves given a square of functors
\begin{equation*}
    \cd{
        \A
            \ar[r]^-H
            \ar[d]_F &
        \C
            \ar[d]^G \\
        \B
            \ar[r]_-K &
        \D\text,
    }
\end{equation*}
with $F$ an $\Ll$-coalgebra and $G$ an $\Rr$-algebra. The $\Ll$-coalgebra
structure on $F$ provides us with a functor $F^\ast \colon \B \to \A$ and a
natural transformation $\eta \colon 1 \Rightarrow FF^\ast$. Thus we can define
a functor $HF^\ast \colon \B \to \C$ and a natural transformation
\begin{equation*}
    \cd{
        \B \ar[rr]^{HF^\ast} \ar[dr]_K & {}\rtwocell{d}{\alpha} & \C \ar[dl]^G \\ & \D\text;
    }
\end{equation*}
indeed, we have $GHF^\ast = KFF^\ast$, and so can take $\alpha = K \eta \colon
K \Rightarrow KFF^\ast$. Now using the $\Rr$-algebra structure on $G$, we may
factorise this 2-cell as:
\begin{equation*}
    \cd[@R+1em]{
        \B \ar@/^0.8em/[rr]^{HF^\ast} \ar@/_0.8em/[rr]_J \ar[dr]_K & {}\utwocell[0]{d}{\overline \alpha} \twoeq{d} & \C \ar[dl]^G \\ & \D\text,
    }
\end{equation*}
where $J$ is given by reindexing $HF^\ast$ along $\alpha$. It is not hard to
see that this functor $J \colon \C \to \D$ is precisely the fill-in specified
by equation \eqref{fillineq} above.
\end{Ex}

\begin{Rk}\label{underlyingplain} It follows from the observations of \S \ref{factdefn} and \S
\ref{worthdefn} that any natural w.f.s.\ $(\Ll, \Rr)$ on a category $\C$ has an
underlying plain w.f.s. For if we define $\ELL$ to be the class of arrows in
$\C$ which admit some $\Ll$-coalgebra structure and $\R$ to be the class of
arrows admitting some $\Rr$-algebra structure, then $(\ELL, \R)$ will satisfy
all the axioms required of a w.f.s., expect possibly for closure under
retracts. So we take $\bar \ELL$ and $\bar \R$ to be the respective
retract-closures of $\ELL$ and $\R$; and now the pair $(\bar \ELL, \bar \R)$
gives a w.f.s.\ on $\C$.
\end{Rk}

\pgph \label{strengthenslightly} It turns to be very useful to strengthen the
notion of natural w.f.s.\ slightly. For this, we consider the natural
transformations $\Pi \colon RR \Rightarrow R$ and $\Sigma \colon L \Rightarrow
LL$ associated to a natural w.f.s.\ $(\Ll, \Rr)$. We may denote their
respective components at $f \in \C^\mathbf 2$ by
\begin{equation*}
    \Pi_f = \cd{
        K\rho_f \ar[d]_{\rho_{\rho_f}} \ar[r]^{\pi_f} & Kf \ar[d]^{\rho_f} \\
        B \ar[r]_{\id_B} & B
    } \qquad \text{and} \qquad
    \Sigma_f = \cd{
        A \ar[d]_{\lambda_f} \ar[r]^{\id_A} & A \ar[d]^{\lambda_{\lambda_f}} \\
        Kf \ar[r]_{\sigma_f} & K\lambda_f\text;
    }
\end{equation*}
again, the arrows written as identities are forced to be so by the (co)monad
axioms. Now, these maps $\sigma_f$ and $\pi_f$ provide us with the components
of a natural transformation $\Delta \colon LR \Rightarrow RL$ whose component
at $f$ is given by:
\begin{equation*}
    \Delta_f = \cd{
        Kf \ar[d]_{\lambda_{\rho_f}} \ar[r]^{\sigma_f} & K\lambda_f \ar[d]^{\rho_{\lambda_f}} \\
        K\rho_f \ar[r]_{\pi_f} & Kf\text.
    }
\end{equation*}
(That this square commutes is a consequence of the (co)monad axioms). We will
say that a natural w.f.s.\ \emph{satisfies the distributivity axiom} if this
natural transformation $\Delta \colon LR \Rightarrow RL$ defines a distributive
law of $\Ll$ over $\Rr$ in the sense of \cite{Beck}. Note that this is a
property of a natural w.f.s., rather than extra structure on it.

\begin{Ex}
We may check that each of the natural w.f.s.'s given so far satisfies the
distributivity axiom.
\end{Ex}

\pgph There are important results about n.w.f.s.'s that are true only if we
restrict to those for which the distributivity axiom holds. Two such results
are Theorem \ref{diamondtwofold} and Theorem \ref{freealgfree} below; and there
is another which allows us to characterise $\Rr$-maps purely in terms of
lifting properties against the $\Ll$-maps, and vice versa. In order that these
results should be valid, we henceforth modify the definition of natural w.f.s.\
to include the requirement that the distributivity axiom should hold.

\section{Free and algebraically-free natural w.f.s.'s}
\pgph Our goal is to use the theory of natural w.f.s.'s to give a categorically
coherent reformulation of the small object argument. As we stated in the
Introduction, this argument provides the means by which, starting from a set of
maps $J$, one may produce a w.f.s.\ \emph{cofibrantly generated} by $J$: that
is, a w.f.s.\ $(\ELL, \R)$ for which $\R = J^\pitchfork$.
\begin{Exs}\label{exs4}
All the weak factorisation systems of Examples \ref{exs1} are cofibrantly
generated:
\begin{itemize}
\item For the w.f.s.\ (Mono, Epi) on $\cat{Set}$, a suitable $J$ is given
    by the set containing the single map $! \colon 0 \to 1$.
\item For the (injective equivalences, isofibrations) w.f.s.\ on
    $\cat{Cat}$, a suitable $J$ is given by the single map $\elt{b} \colon
    1 \to \cat{Iso}$, where $\cat{Iso}$ is the indiscrete category on the
    set $\{a, b\}$.
\item For the w.f.s.\ (anodyne extensions, Kan fibrations) on $\cat{SSet}$,
    a suitable $J$ is given by the set of horn inclusions $\{\Lambda_n^k
    \to \Delta_n\}$.
\end{itemize}
\end{Exs}
To give our reformulation of the small object argument, we will need to provide
a notion of ``cofibrantly generated'' natural w.f.s. However, a careful
analysis reveals two candidates for this notion. In this section, we study
these candidates and their relationship to each other.

\pgph We begin by forming the entities that we have met so far into categories.
Suppose we are given functorial factorisations $F$ and $F' \colon \C^\mathbf 2
\to \C^\mathbf 3$ on $\C$. We define a \emph{morphism of functorial
factorisations} $\alpha \colon F \to F'$ to be a natural transformation $\alpha
\colon F \Rightarrow F'$ which upon whiskering with $d_1 \colon \C^\mathbf 3
\to \C^\mathbf 2$ becomes the identity transformation $\id_{\C^\mathbf 2}
\Rightarrow \id_{\C^\mathbf 2}$. To give such a morphism is to give a family of
maps $\alpha_f \colon Kf \to K'f$, natural in $f$, and making diagrams of the
following form commute:
\[\cd{
 & A \ar[dl]_{\lambda_f} \ar[dr]^{\lambda'_f} \\
 Kf \ar[dr]_{\rho_f} \ar[rr]_{\alpha_f} & & K'f \ar[dl]^{\rho'_f} \\ & B\text.}
\]

Suppose now that $F$ and $F'$ underlie natural w.f.s.'s $(\Ll, \Rr)$ and
$(\Ll', \Rr')$ on $\C$, and consider a morphism of functorial factorisations
$\alpha \colon F \to F'$. By whiskering the natural transformation $\alpha
\colon F \Rightarrow F'$ with the other two face maps $d_0, d_2 \colon
\C^\mathbf 3 \to \C^\mathbf 2$, we induce natural transformations $\alpha_l
\colon L \Rightarrow L'$ and $\alpha_r \colon R \Rightarrow R'$; and we will
say that $\alpha \colon F \to F'$ is a \emph{morphism of natural w.f.s.'s} just
when $\alpha_l$ is a comonad morphism and $\alpha_r$ a monad morphism.

\pgph Let us write $\cat{NWFS}(\C)$ for the category of n.w.f.s.'s on $\C$. We
may define a ``semantics'' functor $\G \colon \cat{NWFS}(\C) \to
\cat{CAT}/\C^\mathbf 2$, which sends a n.w.f.s.\ $(\Ll, \Rr)$ to its category
of $\Ll$-coalgebras $\Lmap$, equipped with the forgetful functor into
$\C^\mathbf 2$; and sends a morphism $\alpha \colon (\Ll, \Rr) \to (\Ll',
\Rr')$ of n.w.f.s.'s to the morphism
\begin{equation*}
\cd{
    \Lmap \ar[rr]^{(\alpha_l)_\ast} \ar[dr]_{U_{\Ll}} & &
    \Ll'\text-\cat{Map} \ar[dl]^{U_{\Ll'}} \\ & \C^\mathbf 2
}
\end{equation*}
of $\cat{CAT}/\C^\mathbf 2$. Here $(\alpha_l)_\ast$ is the functor which sends
an $\Ll$-coalgebra $x \colon X \to LX$ to the $\Ll'$-coalgebra
\begin{equation*}
    X \xrightarrow x LX \xrightarrow{(\alpha_l)_X} L'X\text.
\end{equation*}
Our first candidate for the notion of ``cofibrantly generated'' n.w.f.s.\ is
now:
\begin{Defn}\label{cand1}
Let $I \colon \J \to \C^\mathbf 2$ be an object of $\cat{CAT}/\C^\mathbf 2$,
with $\J$ small; and let $(\Ll, \Rr)$ be a n.w.f.s.\ on $\C$. We will say that
$(\Ll, \Rr)$ is \emph{free on $\J$}\footnote{Here we commit the usual abuse of
notation in denoting a category $I \colon \J \to \C^\mathbf 2$ over $\C^\mathbf
2$ merely by its domain category $\J$.} if we can provide a morphism
\begin{equation*}
\cd[@!C@C-1em]{
    \J \ar[rr]^{\eta} \ar[dr]_{I} & &
    \Ll\text-\cat{Map} \ar[dl]^{U_{\Ll}} \\ & \C^\mathbf 2
}
\end{equation*}
of $\cat{CAT} / \C^\mathbf 2$ which exhibits $(\Ll, \Rr)$ as a reflection of
$I$ along $\G$: which is to say that, for any n.w.f.s.\ $(\Ll', \Rr')$ on $\C$
and functor $F \colon \J \to \Ll'\text-\cat{Map}$ over $\C^\mathbf 2$, there is
a unique morphism of n.w.f.s.'s $\alpha \colon (\Ll, \Rr) \to (\Ll', \Rr')$ for
which $F = (\alpha_l)_\ast \circ \eta$.
\end{Defn}

\begin{Rk}
There is a dual semantics functor $\H \colon \cat{NWFS}(\C) \to
(\cat{CAT}/\C^\mathbf 2)^\op$, which sends a n.w.f.s.\ to its category of
$\Rr$-algebras: and a corresponding notion of an n.w.f.s.\ being \emph{cofree}
on $\J$. However, being cofree is significantly less common than being free,
primarily because the conditions under which we will construct free
n.w.f.s.'s~--~typically, local presentability or local boundedness~--~are much
more prevalent than their duals.
\end{Rk}

\pgph Whilst Definition \ref{cand1} is natural from a categorical perspective,
it has an obvious drawback: it provides no analogue of the equation $\R =
J^\pitchfork$ which a cofibrantly generated w.f.s.\ satisfies. Definition
\ref{cand2}, our second candidate for the notion of ``cofibrantly generated''
n.w.f.s., will rectify this. Before we can give it, we will need a preliminary
result.

\begin{Prop}
Let $\C$ be a category. Then the Galois connection ${}^{\pitchfork}(\ )
\mathbin \dashv (\ )^\pitchfork$ induced by the notion of weak orthogonality
may be lifted to an adjunction
\begin{equation*}
  \cd{\cat{CAT}/\C^\mathbf 2 \ar@<4pt>[r]^-{{}^{\pitchfork}(\thg)} \ar@{}[r]|-{\bot} & (\cat{CAT}/\C^\mathbf 2)^\op \ar@<4pt>[l]^-{(\thg)^\pitchfork}}\text.
\end{equation*}
\end{Prop}

\begin{proof}
First we give the functor $(\thg)^\pitchfork \colon (\cat{CAT}/\C^\mathbf
2)^\op \to \cat{CAT}/\C^\mathbf 2$. On objects, this sends a category $U \colon
\A \to \C^\mathbf 2$ over $\C^\mathbf 2$ to the following category
$\A^\pitchfork$ over $\C^\mathbf 2$. Its objects are pairs $(g, \phi)$, where
$g$ is a morphism of $\C$ and $\phi$ is a coherent choice of lifting against
every element of $\A$: which is to say, a mapping which to each object $a \in
\A$ and square
\begin{equation}\label{fillin2}
    \cd{
        A
            \ar[r]^-h
            \ar[d]_{Ua} &
        C
            \ar[d]^g \\
        B
            \ar[r]_-k &
        D
    }
\end{equation}
in $\C$, assigns a morphism $\phi(a, h, k) \colon B \to C$ making both
triangles commute, and subject to the following naturality condition: if we are
given a morphism $\sigma \colon a \to a'$ of $\A$ whose image under $U$ is the
morphism
\begin{equation*}
    \cd{
        A
            \ar[r]^-s
            \ar[d]_{Ua} &
        A'
            \ar[d]^{Ua'} \\
        B
            \ar[r]_-t &
        B'\text,
    }
\end{equation*}
of $\C^\mathbf 2$, then we have $\phi(a, hs, kt) = \phi(a', h, k) \circ t$. A
morphism of $\A^\pitchfork$ from $(g, \phi)$ to $(g', \phi')$ is a morphism
$(u, v) \colon g \to g'$ of $\C^\mathbf 2$ which respects the choice of
liftings in $\phi$ and $\phi'$, in the sense that the equation $u \circ \phi(a,
h, k) = \phi'(a, uh, vk)$ holds for all suitable $a$, $h$ and $k$. The functor
exhibiting $\A^\pitchfork$ as a category over $\C^\mathbf 2$ is the evident
forgetful functor.

This defines $(\thg)^\pitchfork$ on objects of $\cat{CAT}/\C^\mathbf 2$; and to
extend this definition to morphisms, we consider a further category $\B$ over
$\C^\mathbf 2$ and a functor $F \colon \A \to \B$ over $\C^\mathbf 2$: from
which we obtain a map $F^\pitchfork \colon \B^\pitchfork \to \A^\pitchfork$
over $\C^\mathbf 2$ by sending the object $(c, \phi(\thg, \mathord \ast,
\mathord ?))$ of $\B^\pitchfork$ to the object $(c, \phi(F(\thg), \mathord
\ast, \mathord ?))$ of $\A^\pitchfork$.

We define the functor ${}^\pitchfork(\thg)$ in the same way as
$(\thg)^\pitchfork$, but with $Ua$ and $g$ swapped around in equation
\eqref{fillin2}. It remains only to exhibit the adjointness
${}^\pitchfork(\thg) \dashv (\thg)^\pitchfork$: for which it is easy to see
that, given categories $U \colon \A \to \C^\mathbf 2$ and $V \colon \B \to
\C^\mathbf 2$ over $\C^\mathbf 2$, we may identify both
\begin{equation*}
    \text{functors } \A \to {}^\pitchfork \B \qquad \text{and} \qquad \text{functors } \B \to \A^\pitchfork
\end{equation*}
over $\C^\mathbf 2$ with ``$(\A, \B)$-lifting operations'': that is, functions
$\psi$ which given an object $a \in \A$, an object $b \in \B$ and a commuting
square
\begin{equation*}
    \cd{
        A
            \ar[r]^-h
            \ar[d]_{Ua} &
        C
            \ar[d]^{Vb} \\
        B
            \ar[r]_-k &
        D\text,
    }
\end{equation*}
provide a morphism $\psi(a, b, h, k) \colon B \to C$ making both triangles
commute; and subject to the obvious naturality condition with respect to
morphisms of both $\A$ and $\B$.
\end{proof}

In particular, we see from \S \ref{worthdefn} that any any n.w.f.s.\ $(\Ll,
\Rr)$ comes equipped with a privileged $(\Lmap$, $ \Rmap)$-lifting operation:
which by the above proof, we may view as a privileged morphism $\textsf{lift}
\colon \Rmap \to \Lmap^\pitchfork$ over $\C^\mathbf 2$.

\begin{Defn}\label{cand2}
Let $I \colon \J \to \C^\mathbf 2$ be a category over $\C^\mathbf 2$, and
$(\Ll, \Rr)$ a n.w.f.s.\ on $\C$. We say that $(\Ll, \Rr)$ is
\emph{algebraically-free} on $\J$ if we can provide a morphism $\eta \colon \J
\to \Lmap$ over $\C$ for which the functor
\begin{equation}\label{algfreefunctor}
    \Rmap \xrightarrow{\textsf{lift}} \Lmap^\pitchfork \xrightarrow{\eta^\pitchfork}
    \J^\pitchfork
\end{equation}
is an isomorphism of categories.
\end{Defn}

\begin{Rk}
The terminology we have chosen deliberately recalls the distinction which is
made in \cite{Ke80} between the \emph{free} and the \emph{algebraically-free}
monad generated by a pointed endofunctor. We will partially justify this in
Section \ref{Sec:cfaf}, by showing that algebraic-freeness in our sense can be
seen as a special case of algebraic-freeness in the sense of \cite{Ke80}; and
in the Appendix, where we prove the implication ``algebraically-free
$\Rightarrow$ free'' for n.w.f.s.'s.

However, there are some results of \cite{Ke80} which the author has been unable
to find an analogue of: in particular,  he has been unable to produce either
positive or negative results about the implication ``free $\Rightarrow$
algebraically-free''. The corresponding implication does not hold in the theory
of monads; and whilst it seems unlikely that it should hold here either, a
proof of this fact has been elusive. Despite this, we will be able to show in
Section \ref{Sec:cfaf} that any free n.w.f.s.\ which we come across \emph{in
mathematical practice} will be algebraically-free.
\end{Rk}

\begin{Exs}\label{exs5}
The natural w.f.s.\ on $\cat{Set}$ which factors $f \colon X \to Y$ through $X
+ Y$ is algebraically-free: we let $\J$ be the terminal category and let $I
\colon \J \to \cat{Set}^\mathbf 2$ pick out the object $! \colon 0 \to 1$. It
is now easy to see that the category $\J^\pitchfork$ consists precisely of the
$\Rr$-algebras: morphisms $g \colon C \to D$ equipped with a splitting $g^\ast
\colon D \to C$.

However, none of the other natural w.f.s.'s described in Examples \ref{exs2}
are free or algebraically-free: and this despite being close relatives of plain
w.f.s.'s which are cofibrantly generated. The problem for these examples is
that, although an $\Rr$-map structure affirms the existence of certain
liftings, it also asserts certain coherence conditions between those liftings,
which cannot be expressed in the language of orthogonality.

A fair intuition is that the (algebraically)-free natural w.f.s.'s are the
natural w.f.s.'s which may be specified by a ``signature'' $\J$ of lifting
properties; but subject to no ``equations'' between these liftings.
\end{Exs}

We may relate the notion of algebraically-free n.w.f.s.\ quite directly to that
of cofibrantly generated w.f.s., if we assume the axiom of choice in our
metatheory:

\begin{Propstar}\label{relateback}
Let $\C$ be a category and $J$ a set of maps in $\C$; and let $\J$ denote the
set $J$, viewed as a discrete subcategory of $\C^\mathbf 2$. If the
algebraically-free n.w.f.s.\ $(\Ll, \Rr)$ on $\J \hookrightarrow \C^\mathbf 2$
exists, then its underlying plain w.f.s.\ $(\bar \ELL, \bar \R)$ is the w.f.s.\
cofibrantly generated by $J$.
\end{Propstar}
\begin{proof}
Recall from \S \ref{underlyingplain} that the class of maps $\R$ consists of
those maps in $\C$ admitting some $\Rr$-algebra structure; and that $\bar \R$
consists of all retracts of maps in $\R$. We are required to show that $\bar \R
= J^\pitchfork$; and since $J^\pitchfork$ is easily seen to be closed under
retracts, it will suffice to show that $\R = J^\pitchfork$.

Now, since $(\Ll, \Rr)$ is algebraically-free on $\J$, we have $\Rmap \cong
\J^\pitchfork$ over $\C^\mathbf 2$; and so a morphism $f \in \C^\mathbf 2$ will
admit an $\Rr$-algebra structure, and thus lie in $\R$, just when it can be
lifted through the forgetful functor $\J^\pitchfork \to \C^\mathbf 2$. But an
object of $\J^\pitchfork$ consists of a map of $\C$ equipped with a choice of
lifting against every map in the set $J$, subject to no further coherence
conditions; and so, if we allow ourselves the axiom of choice, $f$ will admit a
lifting through $\J^\pitchfork$ just when $f \in J^\pitchfork$. Thus we have
that $\R = J^\pitchfork$ as desired.
\end{proof}

\section{Constructing free natural w.f.s.'s}\label{construct}
\pgph We now ready to give our analogue of the small object argument, which
will be a general apparatus by means of which we can construct free, and even
algebraically-free, n.w.f.s.'s on a category $\C$.

For our argument to work, we will at least require $\C$ to be cocomplete: but
in order to guarantee the convergence of certain transfinite sequences we
construct, we must impose some further ``smallness'' property on $\C$.

\pgph Given a regular cardinal $\alpha$, we say that $X \in \C$ is
$\alpha$-\emph{presentable} if the representable functor $\C(X, \thg) \colon \C
\to \cat{Set}$ preserves $\alpha$-filtered colimits. The first smallness
property we may consider on $\C$ is that:
\begin{quote}
 (*) For every $X \in \C$, there is an $\alpha_X$ for which $X$ is
$\alpha_X$-presentable.
\end{quote}
\noindent This is certainly the case for any category $\C$ which is
\emph{locally presentable} in the sense of \cite{lpk}. However, it does not
obtain in categories such as the category of topological spaces, the category
of Hausdorff topological spaces, or the category of topological groups: and
since we would like our argument to be valid in such contexts, we will require
a more general notion of smallness.

\pgph Recall that a strong factorisation system $(\E, \M)$ on $\C$ is said to
be \emph{proper} if every $\E$-map is an epimorphism and every $\M$-map a
monomorphism; and is said to be \emph{well-copowered} if every object of $\C$
possesses, up-to-isomorphism, a mere set of $\E$-quotients. We say that an
object $X$ is $\alpha$-\emph{bounded} with respect to a proper $(\E, \M)$ if
$\C(X, \thg)$ preserves $\alpha$-filtered unions of $\M$-subobjects (in the
sense of sending them to $\alpha$-filtered unions of sets). The second
smallness property we consider on $\C$ supposes some proper, well-copowered
$(\E, \M)$, and says that:
\begin{quote}
(\dag) For every $X \in \C$, there is an $\alpha_X$ for which $X$ is
$\alpha_X$-bounded with respect to $(\E, \M)$.
\end{quote}

\noindent $\cat{Top}$, $\cat{Haus}$ and $\cat{TopGrp}$ all satisfy (\dag), with
$\M $ = the subspace inclusions in the first two cases, and $\M = $ the
inclusion of subgroups which are also subspaces in the third.

\vskip\baselineskip

We may now state the main result of the paper.

\begin{Thm}\label{mainthm}
Let $\C$ be a cocomplete category satisfying either (*) or (\dag), and let $I
\colon \J \to \C^\mathbf 2$ be a category over $\C^\mathbf 2$ with $\J$ small.
Then the free n.w.f.s.\ on $\J$ exists, and is algebraically-free on $\J$.
\end{Thm}

In this section, we will prove freeness: in the next, algebraic-freeness.

\pgph We begin by factorising the semantics functor $\G$ through a pair of
intermediate categories. The first is the category $\cat{LNWFS}(\C)$ of ``left
halves of n.w.f.s.'s''. Its objects $(F, \Ll)$ are functorial factorisations
$F$ on $\C$ together with an extension of the corresponding $(L, \Phi)$ to a
comonad $\Ll$; and its morphisms are maps of functorial factorisations which
respect the comonad structure. There is an obvious functor $\G_1 \colon
\cat{NWFS}(\C) \to \cat{LNWFS}(\C)$ sending $(\Ll, \Rr)$ to $(F, \Ll)$.

The second category we consider is $\cat{Cmd}(\C^\mathbf 2)$, the category of
comonads on $\C^\mathbf 2$. We have a functor $\G_2 \colon \cat{LNWFS}(\C) \to
\cat{Cmd}(\C^\mathbf 2)$, which sends $(F, \Ll)$ to $\Ll$; and we have the
semantics functor $\G_3 \colon \cat{Cmd}(\C^\mathbf 2) \to \cat{CAT} /
\C^\mathbf 2$ which sends a comonad to its category of coalgebras, and a
comonad morphism $\gamma \colon \mathsf C \to \mathsf C'$ to $\gamma_\ast
\colon \mathsf C\text-\cat{Coalg} \to \mathsf C'\text-\cat{Coalg}$. We now have
that
\begin{equation*}
    \G \quad = \quad \cat{NWFS}(\C) \xrightarrow{\G_1} \cat{LNWFS}(\C) \xrightarrow{\G_2}
    \cat{Cmd}(\C^\mathbf 2) \xrightarrow{\G_3} \cat{CAT} / \C^\mathbf 2\text,
\end{equation*}
so that we may give a reflection along $\G$ by giving a reflection along each
functor $\G_1$, $\G_2$ and $\G_3$ in turn. For $\G_3$, we have the following
well-known result, which was first stated at this level of generality by Dubuc
\cite{Dbc}; but see also \cite{triples}.

\begin{Prop}\label{Kanext}
Let $\C$ be cocomplete, and let $U \colon \A \to \C^\mathbf 2$ be a small
category over $\C^\mathbf 2$. Then $\A$ admits a reflection along $\G_3 \colon
\cat{Cmd}(\C^\mathbf 2) \to \cat{CAT} / \C^\mathbf 2$.
\end{Prop}
\begin{proof}
Because $\A$ is small and $\C^\mathbf 2$ cocomplete (since $\C$ is), we can
form the left Kan extension of $U$ along itself:
\begin{equation*}
    \cd{
    & \A \ar[dr]^U \ar[dl]_U \ltwocell{d}{\theta} \\
  \C^\mathbf 2 \ar[rr]_{\Lan_U(U)} & & \C^\mathbf 2\text,
}
\end{equation*}
whose defining property is that $\theta$ should provide the unit for a
representation
\begin{equation*}
    [\C^\mathbf 2, \C^\mathbf 2](\Lan_U(U), \thg) \cong [\A, \C^\mathbf 2](U, (\thg) \circ U)\text.
\end{equation*}
In particular, corresponding to the identity transformation $\id_U \colon U
\Rightarrow U$, we have a natural transformation $\epsilon \colon \Lan_U(U)
\Rightarrow \id_{\C^\mathbf 2}$; whilst corresponding to the composite natural
transformation
\begin{equation*}
    U \xrightarrow{\theta} \Lan_U(U) \circ U \xrightarrow{\Lan_U(U) \circ \theta} \Lan_U(U) \circ
    \Lan_U(U) \circ U
\end{equation*}
we have a natural transformation $\Delta \colon \Lan_U(U) \Rightarrow \Lan_U(U)
\circ \Lan_U(U)$. It is now easy to check that $\epsilon$ and $\Delta$ make
$\Lan_U(U)$ into a comonad on $\C^\mathbf 2$, the so-called \emph{density
comonad} of $U$. This has the property that comonad morphisms $(\Lan_U(U),
\epsilon, \Delta) \to \mathsf T$ are in bijection with left coactions of
$\mathsf T$ on $U$, which in turn are in bijection with liftings of $U \colon
\A \to \C^\mathbf 2$ through the category of $\mathsf T$-coalgebras: and this
is precisely the universal property for $\Lan_U(U)$ to be a reflection of $U$
along $\G_3$.
\end{proof}

Next, we consider reflections along $\G_2 \colon \cat{LNWFS}(\C) \to
\cat{Cmd}(\C^\mathbf 2)$. These exist under very mild hypotheses indeed:

\begin{Prop}\label{reflection}
If $\C$ has pushouts, then $\G_2 \colon \cat{LNWFS}(\C) \to
\cat{Cmd}(\C^\mathbf 2)$ has a left adjoint.
\end{Prop}
\begin{proof}
Let us say that an endofunctor $F \colon \C^\mathbf 2 \to \C^\mathbf 2$
\emph{preserves domains} if $\dom \circ F = \dom$. Given two such endofunctors
$F$ and $F'$, we will say that a natural transformation $\alpha$ between them
\emph{preserves domains} if $\dom \circ \alpha = \id_\dom$. Finally, we will
say that a comonad $(T, \epsilon, \Delta)$ on $\C^\mathbf 2$ \emph{preserves
domains} if $T$, $\epsilon$ and $\Delta$ all preserve domains.

It is now a simple but instructive exercise to show that $\cat{LNWFS}(\C)$ is
isomorphic to the full subcategory of $\cat{Cmd}(\C^\mathbf2)$ whose objects
are the domain-preserving comonads. Thus the Proposition will follow if we can
show this subcategory to be reflective.

To do this, we first observe that there is a strong factorisation system on
$\C^\mathbf 2$ whose left class $\P$ consists of the pushout squares, and whose
right class consists $\D$ of the squares whose domain component is an
isomorphism. In fact, if we make a choice of pushouts in $\C$, then we obtain a
functorial factorisation of every map into a pushout square followed by a
square whose domain component is an \emph{identity}.

We can lift the factorisation system $(\P, \D)$ to one of the same name on
$[\C^\mathbf 2, \C^\mathbf 2]$; and the accompanying functorial factorisation
lifts too, allowing us to factor every map of $[\C^\mathbf 2, \C^\mathbf 2]$ as
a map whose components are pushouts, followed by one whose domain components
are identities.

Suppose now that we are given a comonad $\mathsf S = (S, \epsilon, \Delta)$ on
$\C^\mathbf 2$: we construct its reflection into domain-preserving comonads as
follows. We start by factorising the counit of $\mathsf S$ as
\begin{equation*}
    \epsilon = S \xRightarrow{\phi} \hat S \xRightarrow{\hat \epsilon} \id_{\C^\mathbf
    2}\text,
\end{equation*}
where the components of $\phi$ are pushout squares, and the domain components
of $\hat \epsilon$ are identities. From this latter fact, we deduce that both
$\hat S$ and $\hat \epsilon$ preserve domains. We now consider the following
diagram:
\begin{equation*}
\cd{
    S \ar@2[r]^{\Delta} \ar@2[d]_{\phi} &
    SS \ar@2[r]^{\phi \phi} &
    \hat S \hat S \ar@2[d]^{\hat \epsilon \hat S} \\
    \hat S \ar@2[rr]_{\id_{\hat{S}}} & & \hat S\text.
}
\end{equation*}
Since $\phi$ is in $\P$, and $\hat \epsilon \hat S$ in $\D$, we obtain by
orthogonality a unique diagonal fill-in $\hat \Delta \colon \hat S \Rightarrow
\hat S \hat S$. Since both $\hat \epsilon \hat S$ and $\id_{\hat S}$ are
domain-preserving, we deduce that $\hat \Delta$ is too.

A little calculus with the unique diagonal fill-in property and the comonad
axioms for $(S, \epsilon, \Delta)$ now yields the comonad axioms for
$\hat{\mathsf S} = (\hat S, \hat \epsilon, \hat \Delta)$; and it is immediate
that $\phi \colon S \Rightarrow \hat S$ then satisfies the necessary axioms for
it to lift to a comonad morphism $\phi \colon \mathsf S \to \hat{\mathsf S}$.

We claim that this $\phi$ provides the desired reflection of $\mathsf S$ into
domain-pre\-ser\-ving comonads. Indeed, suppose we are given another
domain-preserving comonad $\mathsf T = (T, e, D)$, and a morphism of comonads
$\psi \colon \mathsf S \to \mathsf T$. Then we have the following commutative
square:
\begin{equation*}
\cd{
    S \ar@2[r]^{\psi} \ar@2[d]_{\phi} &
    T \ar@2[d]^{e} \\ \hat S \ar@2[r]_{\hat \epsilon} & \id_{\C^\mathbf 2}\text.
}
\end{equation*}
The map $\phi$ is in $\P$, and $e$ is in $\D$: so by orthogonality, we induce a
unique natural transformation $\hat \psi \colon \hat S \Rightarrow T$. The
comonad morphism axioms for $\hat \psi$ now follow from the axioms for $\psi$
and uniqueness of diagonal fill-ins.
\end{proof}

\pgph We have thus reduced the problem of constructing free n.w.f.s.'s to the
problem of constructing reflections along $\G_1 \colon \cat{NWFS}(\C) \to
\cat{LNWFS}(\C)$. The key to constructing these will be to exhibit a monoidal
structure on $\cat{LNWFS}(\C)$ whose corresponding category of monoids is
isomorphic to $\cat{NWFS}(\C)$.

We will deduce the existence of this monoidal structure from a more general
result characterising natural w.f.s.'s on $\C$ as \emph{bialgebra} objects in
the category  of functorial factorisations on $\C$. Now, usually when one
considers bialgebra objects in a category, it is with reference to a symmetric
or braided monoidal structure on that category: but here we will need something
slightly more general.

\pgph By a \emph{two-fold monoidal category} \cite{iter1}, we mean a category
$\V$ equipped with two monoidal structures $(\otimes, I, \alpha, \lambda,
\rho)$ and $(\odot, \bot, \alpha', \lambda', \rho')$ in such a way that the
functors $\odot \colon \V \times \V \to \V$ and $\bot \colon 1 \to \V$,
together with the natural transformations $\alpha'$, $\lambda'$ and $\rho'$,
are lax monoidal with respect to the $(\otimes, I)$ monoidal structure.

Of course, being lax monoidal is not merely a property of a functor, but extra
structure on it: and in this case, the extra structure amounts to giving maps
\begin{equation*}
    m \colon \bot \otimes \bot \to \bot\text, \quad c \colon I \to I \odot I \quad \text{and} \quad j \colon I \to \bot
\end{equation*}
making $(\bot, j, m)$ into a $\otimes$-monoid and $(I, j, c)$ into a
$\odot$-comonoid; together with a natural family of maps
\[z_{A,B,C,D} \colon (A \odot B) \otimes (C \odot D) \to (A \otimes C) \odot (B \otimes D)\]
obeying six coherence laws. It follows that $\otimes$ and $I$ are oplax
monoidal with respect to the $(\odot, \bot)$ monoidal structure; and in fact,
we may take this as an alternative definition of two-fold monoidal category.

\begin{Exs}
\begin{itemize*}
\item Any braided or symmetric monoidal category is two-fold monoidal, with the
    two monoidal structures coinciding; the maps $z_{A, B, C, D}$ are built
    from braidings/symmetries and associativity isomorphisms: c.f.\ \cite{Braided}.
\item If $\V$ is a cocomplete symmetric monoidal category, then the functor
    category $[X \times X, \V]$ has a two-fold monoidal structure. The first
    monoidal structure $(\otimes, I)$ is given by matrix multiplication, whilst
    the second structure $(\odot, \bot)$ is given pointwise.
\item Similarly, if $\V$ is a cocomplete symmetric monoidal category, then the
    functor category $[\mathbb N, \V]$ has a two-fold monoidal structure on it.
    The first monoidal structure $(\otimes, I)$ is the \emph{substitution tensor
    product}, with unit given by $I(n) = 0$ for $n \neq 1$ and $I(1) = I$; and binary tensor given by
    \begin{equation*}
        (F \otimes G)(n) = \sum_{\substack{m, k_1, \dots, k_m \\ k_1 + \dots + k_m = n}} F(m) \otimes G(k_1) \otimes \dots \otimes G(k_m)\text.
    \end{equation*}
    The second monoidal structure $(\odot, \bot)$ is again given pointwise.
\end{itemize*}
\end{Exs}

Further examples and applications to topology may be found in \cite{iter1,
iter2}.

\pgph A two-fold monoidal category $(\V, \otimes, I, \odot, \bot)$ provides a
suitable environment to define a notion of bialgebra. Indeed, because the
$\odot$-monoidal structure is lax monoidal with respect to the
$\otimes$-structure, it lifts to the category $\cat{Mon}_\otimes(\V)$ of
$\otimes$-monoids in $\V$. Thus we define the category of \emph{$(\otimes,
\odot)$-bialgebras} to be
\begin{equation*}
    \cat{Bialg}_{\otimes, \odot}(\V) := \cat{Comon}_\odot(\cat{Mon}_\otimes(\V))\text.
\end{equation*}

Now, because the $\otimes$-monoidal structure is also oplax monoidal with
respect to the $\odot$-monoidal structure, it lifts to the category of
$\odot$-comonoids in $\V$; and thus we obtain an alternative definition of
bialgebra by setting
\begin{equation*}
    \cat{Bialg}'_{\otimes, \odot}(\V) := \cat{Mon}_\otimes(\cat{Comon}_\odot(\V))\text.
\end{equation*}

However, it is not hard to see that these two constructions yield isomorphic
results. Indeed, to give an object of either $\cat{Bialg}(\V)$ or
$\cat{Bialg}'(\V)$ is to give an object $A$ of $\V$; maps $\eta \colon I \to A$
and $\mu \colon A \otimes A \to A$ making it into a $\otimes$-monoid; and maps
$\epsilon \colon A \to \bot$ and $\delta \colon A \to A \odot A$ making it into
a $\odot$-comonoid; all subject to the commutativity of the following four
diagrams:
\begin{equation}
\label{bialg}
\begin{gathered}
\cd{
 I \ar[r]^-\eta \ar[d]_{c} &
 A \ar[d]^\Delta \\
 I \odot I \ar[r]_-{\eta \odot \eta} &
 A \odot A\text,
} \qquad \cd{
 A \otimes A \ar[r]^-\mu \ar[d]_{\epsilon \otimes \epsilon} &
 A \ar[d]^\epsilon \\
 \bot \otimes \bot \ar[r]_-m & \bot\text,
} \qquad \cd{
 & A \ar[dr]^\epsilon \\
 I \ar[ur]^\eta \ar[rr]_j & & \bot\text,
}
\\
\cd[@C+1.5em]{
 A \otimes A
   \ar[rr]^\mu \ar[d]_{\Delta \otimes \Delta} & &
 A \ar[d]^\Delta \\
 (A \odot A) \otimes (A \odot A)
 \ar[r]_-{z_{A, A, A, A}} &  (A \otimes A) \odot (A \otimes A)
 \ar[r]_-{\mu \odot \mu} &
 A \odot A\text.
 }
\end{gathered}
\end{equation}
Likewise, to give a morphism of either $\cat{Bialg}(\V)$ or $\cat{Bialg}'(\V)$
is to give a map $f \colon A \to B$ of $\V$ which is both a monoid morphism and
a comonoid morphism. We may summarise this by saying that, in the following
diamond of forgetful functors
\begin{equation}\label{diamond}
    \cd{
    & \cat{Bialg}_{\otimes, \odot}(\V) \ar[dl] \ar[dr] \\
    \cat{Comon}_\odot(\V) \ar[dr] & &
    \cat{Mon}_\otimes(\V)\text, \ar[dl] \\ &
    \V
    }
\end{equation}
each west-pointing arrow forgets monoid structure, and each east-pointing arrow
forgets comonoid structure.

\begin{Exs}
\begin{itemize*}
\item If view a braided or symmetric monoidal category as a two-fold
    monoidal category, then a bialgebra in our sense is precisely a bialgebra
    in the usual sense.
\item In the two-fold monoidal category $[X \times X, \V]$, a $\otimes$-monoid
    is a $\V$-category with object set $X$; an $\odot$-comonoid is an $X \times
    X$-indexed family of comonoids in $\V$; and a $(\otimes, \odot)$-bialgebra
    is a \emph{comonoidal} $\V$-category with object set $X$: which we may view either
    as a comonoid in $\V\text-\cat{Cat}$, or as a $\V$-category whose homsets are
    comonoids and whose unit and composition maps are comonoid morphisms.
\item A bialgebra in the two-fold monoidal category $[\mathbb N, \V]$
    is what is sometimes called a \emph{Hopf operad}: namely, an operad in $\V$
    whose objects of $n$-ary operations are comonoids; and whose substitution
    maps are morphisms of comonoids.
\end{itemize*}
\end{Exs}

Bialgebras in two-fold monoidal categories play a central role in
    recent work \cite{lamarche} of Fran\c cois Lamarche.

\pgph Let us write $\cat{FF}(\C)$ for the category of functorial factorisations
on $\C$, and let us write $\cat{RNWFS}(\C)$ for the category dual to
$\cat{LNWFS}(\C)$: so its objects are pairs $(F, \Rr)$ of a functorial
factorisation $F$ on $\C$ together with an extension of the corresponding $(R,
\Lambda)$ to a monad.

\begin{Thm}\label{diamondtwofold}
There is a two-fold monoidal structure on $\cat{FF}(\C)$ such that the diamond
of forgetful functors \eqref{diamond} is, up-to-isomorphism, the diamond of
forgetful functors
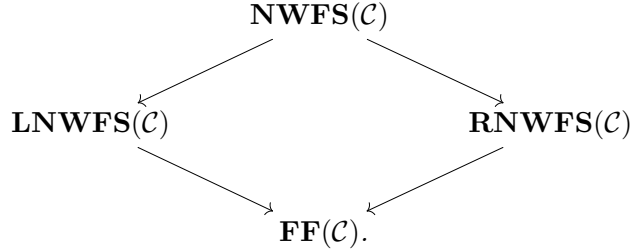
\begin{equation*}
    \cd{
    & \cat{NWFS}(\C)\ar[dl] \ar[dr] \\
    \cat{LNWFS}(\C) \ar[dr] & &
    \cat{RNWFS}(\C) \ar[dl] \\ &
    \cat{FF}(\C)\text{\rlap.}
    }
\end{equation*}
\end{Thm}
\begin{proof}
We begin by exhibiting two strict monoidal structures on $\cat{FF}(\C)$. We do
this by describing two different categories which are both isomorphic to
$\cat{FF}(\C)$, and which both admit obvious strict monoidal structures: then
by transport of structure, we induce the required monoidal structures on
$\cat{FF}(\C)$.

The first category we consider is the category of domain-preserving copointed
endofunctors and copointed endofunctor maps on $\C^\mathbf 2$. It is easy to
see that this category is isomorphic to $\cat{FF}(\C)$; and that it has a
strict monoidal structure $(\odot, \bot)$ on it, with unit
\begin{equation*}
    \bot = (\id_{\id_{\C^\mathbf 2}} \colon \id_{\C^\mathbf 2} \Rightarrow \id_{\C^\mathbf 2})
\end{equation*}
and tensor product
\begin{equation*}
    (\Phi \colon L \Rightarrow \id_{\C^\mathbf 2}) \odot (\Phi' \colon L' \Rightarrow \id_{\C^\mathbf 2}) = (\Phi \Phi' \colon LL' \Rightarrow \id_{\C^\mathbf 2})\text.
\end{equation*}

When we transport this along the isomorphism with $\cat{FF}(\C)$, we obtain the
following monoidal structure. The unit $\bot$ is the functorial factorisation
\begin{equation*}
    X \xrightarrow f Y \quad \mapsto \quad X \xrightarrow f Y \xrightarrow{\id_Y} Y
\end{equation*}
and the tensor product $F' \odot F$ of two functorial factorisations $F, F'
\colon \C^\mathbf 2 \to \C^\mathbf 3$ is given by
\begin{equation*}
    X \xrightarrow f Y \quad \mapsto \quad X \xrightarrow{\lambda'_{\lambda_f}}  K' \lambda_f \xrightarrow {\rho_f \circ \rho'_{\lambda_f}} Y\text.
\end{equation*}

Furthermore, to give a $\odot$-comonoid structure on some $F \in \cat{FF}(\C)$
is to give a comonoid structure on the corresponding copointed $(L, \Phi)$; but
this is precisely to extend it to a comonad on $\C^\mathbf 2$. Thus we may
identify $\cat{Comon}_\odot(\cat{FF}(\C))$ with $\cat{LNWFS}(\C)$.

The second category we consider is the category of codomain-preserving pointed
endofunctors on $\C^\mathbf 2$. Again, this is isomorphic to $\cat{FF}(\C)$,
and again, it has a strict monoidal structure given by composition. When we
transport this back to $\cat{FF}(\C)$, we obtain the monoidal structure whose
unit $I$ is the functorial factorisation
\begin{equation*}
    X \xrightarrow f Y \quad \mapsto \quad X \xrightarrow{\id_X} X \xrightarrow f Y\text;
\end{equation*}
and whose tensor product $F' \otimes F$ is the functorial factorisation
\begin{equation*}
    X \xrightarrow f Y \quad \mapsto \quad X \xrightarrow{\lambda'_{\rho_f} \circ \lambda_f}  K' \rho_f \xrightarrow {\rho'_{\rho_f}} Y\text.
\end{equation*}

To make $F \in \cat{FF}(\C)$ into a monoid with respect to this monoidal
structure is now to give an extension of the corresponding $(R, \Lambda)$ to a
monad; and so we have $\cat{Mon}_\otimes(\cat{FF}(\C)) \cong \cat{RNWFS}(\C)$
as required.

We next show that these two monoidal structures on $\cat{FF}(\C)$ can be made
into a two-fold monoidal structure. Since $I$ is initial and $\bot$ terminal in
$\cat{FF}(\C)$, for this we need only give the family of interchange maps
$z_{A, B, C, D} \colon (A \odot B) \otimes (C \odot D) \to (A \otimes C) \odot
(B \otimes D)$: and this we do explicitly. The factorisation $(A \odot B)
\otimes (C \odot D)$ sends a map $f \colon X \to Y$ to
\[
X \xrightarrow{\textstyle\lambda^A(\lambda^B(\rho^{C \odot D}_f)) \circ
\lambda^C(\lambda^D_f)} {K}^A(\lambda^B(\rho^{C \odot D}_f))
\xrightarrow{\textstyle\rho^B(\rho^{C \odot D}_f) \circ \rho^A(\lambda^B(\rho^{C \odot D}_f))} Y\text,
\]
where $\rho^{C \odot D}_f$ abbreviates the map $\rho^D_f \circ
\rho^C(\lambda^D_f)$; whilst $(A \otimes C) \odot (B \otimes D)$ sends $f$ to
\[
X \xrightarrow{\textstyle\lambda^A(\rho^C(\lambda^{B \otimes D}_f)) \circ
\lambda^C(\lambda^{B \otimes D}_f)} K^A(\rho^C(\lambda^{B \otimes D}_f))
\xrightarrow{\textstyle\rho^B(\rho^D_f) \circ \rho^A(\rho^C(\lambda^{B \otimes D}_f))} Y\text,
\]
where $\lambda^{B \otimes D}_f$ abbreviates the map $\lambda^B(\rho^D_f) \circ
\lambda^D_f$. To give $z_{A, B, C, D}$ we must therefore give suitable maps
${K}^A(\lambda^B(\rho^{C \odot D}_f)) \to K^A(\rho^C(\lambda^{B \otimes
D}_f))$. For this, we consider the following square:
\[
\cd[@+1.5em@C+2em]{
    K^C(\lambda^D_f)
        \ar[d]_{\lambda^B(\rho^{C \odot D}_f)}
        \ar[r]^-{K^C(\id_X, \lambda^B(\rho^D_f))} &
    K^C(\lambda^{B \otimes D}_f)
        \ar[d]^{\rho^C(\lambda^{B \otimes D}_f)} \\
    K^B(\rho^{C \odot D}_f)
        \ar[r]_-{K^B(\rho^C(\lambda^D_f), \id_Y)} &
    K^B(\rho^D_f)\text.
 }
\]
This square commutes, with both sides equal to
\begin{equation*}
    K^C(\lambda^D_f) \xrightarrow{\rho^C(\lambda^D_f)} K^D f \xrightarrow{\lambda^B(\rho^D_f)} K^B(\rho^D_f)\text,
\end{equation*}
and so we may view it as a morphism $\lambda^B(\rho^{C \odot D}_f) \to
\rho^C(\lambda^{B \otimes D}_f)$ in $\C^\mathbf 2$: applying $K^A$ to which
yields the required map ${K}^A(\lambda^B(\rho^{C \odot D}_f)) \to
K^A(\rho^C(\lambda^{B \otimes D}_f))$ in $\C$. The (extensive) remaining
details are left to the reader.

Thus we have a two-fold monoidal structure $(\otimes, \odot)$ on
$\cat{FF}(\C)$: and to complete the proof, we must show that the corresponding
bialgebras are precisely n.w.f.s.'s on $\C$. But to equip a functorial
factorisation with both a $\otimes$-monoid and an $\odot$-comonoid structure is
to give extensions of the corresponding $(R, \Lambda)$ to a monad $\Rr$, and
the corresponding $(L, \Phi)$ to a comonad $\Ll$; and it is now a short
calculation to show that the bialgebra axioms \eqref{bialg} will hold just when
the distributivity axiom holds for $(\Ll, \Rr)$.
\end{proof}

\pgph This Theorem implies that an object $X \in \cat{LNWFS}(\C)$ will admit a
reflection along the functor $\G_1 \colon \cat{NWFS}(\C) \to \cat{LNWFS}(\C)$
just when the free $\otimes$-monoid on $X$ exists. But since the unit $I$ of
the monoidal structure on $\cat{LWNFS}(\C)$ is also an initial object, to
construct the free monoid on $X$ is equally well to construct the free monoid
on the pointed object $! \colon I \to X$. In order to do this, we may employ a
standard transfinite construction: which we now describe.

\pgph If $\cat{On}$ denotes the category of all small ordinals, then a
\emph{transfinite sequence} in a category $\V$ is a functor $X \colon \cat{On}
\to \V$, whose value at an ordinal $\alpha$ we denote by $X_\alpha$, and whose
value at the unique morphism $\alpha \to \beta$ (for $\alpha \leqslant \beta$)
we denote by $X_{\alpha, \beta} \colon X_\alpha \to X_\beta$. We say that a
transfinite sequence $X \colon \cat{On} \to \V$ \emph{converges} at an ordinal
$\gamma$ if the maps $X_{\alpha, \beta}$ are isomorphisms for all $\beta
\geqslant \alpha \geqslant \gamma$.

Let $\V$ now be a cocomplete monoidal category. Given a pointed object $t
\colon I \to T$ in $\V$, we may form a transfinite sequence $X \colon \cat{On}
\to \V$ which we call the \emph{free monoid sequence} for $(T, t)$. We build
this sequence, together with a family of maps $\sigma_\alpha \colon T \otimes
X_\alpha \to X_{\alpha^+}$, by the following transfinite induction:
\begin{itemize}
\item $X_0 = I$, $X_1 = T$, $X_{0, 1} = t$, and $\sigma_0 = \rho_T \colon T
    \otimes I \to T$;
\item For a successor ordinal $\beta = \alpha^+$, we give $X_{\beta}$ and
    $\sigma_{\beta} \colon T \otimes X_\beta \to X_{\beta^+}$ by the
    following coequaliser diagram:
\[\cd[@R-1.5em]{ &
 X_{\beta}
  \ar@/^0.5em/[dr]^-{t \otimes X_{\beta}}
  \\
 T \otimes X_\alpha
  \ar@/^0.5em/[ur]^-{\sigma_\alpha}
  \ar@/_0.5em/[dr]_-{T \otimes t \otimes X_\alpha} & &
 T \otimes X_{\beta}
  \ar[r]^-{\sigma_{\beta}} &
 X_{\beta^{+}}\text,
  \\ &
 T \otimes T \otimes X_{\alpha}
  \ar@/_0.5em/[ur]_-{T \otimes \sigma_{\alpha}} }\] and give $X_{\beta, \beta^+}$ by the composite $\sigma_\beta \circ (t \otimes X_\beta)$;
\item For a non-zero limit ordinal $\gamma$, we give $X_\gamma$ by
    $\colim_{\alpha < \gamma} X_\alpha$, with connecting maps $X_{\alpha,
    \gamma}$ given by the injections into the colimit. We give
    $X_{\gamma^+}$ and $\sigma_\gamma$ by the following coequaliser
    diagram:
\[
\cd[@C-2em]{
 &  \colim X_{\alpha^+} = X_\gamma
  \ar[dr]^{t \otimes X_\gamma} \\
 \colim (T \otimes X_{\alpha})
  \ar[ur]^{\colim \sigma_\alpha\ }
  \ar[rr]_{\textsf{can}} & &
 T \otimes \colim X_{\alpha} = T \otimes X_{\gamma}
  \ar[rrr]_-{\sigma_{\gamma}} & & &
 X_{\gamma^+}
 }
\]
where ``\textsf{can}'' is the map induced by the cocone $T \otimes X_\alpha
\to T \otimes \colim X_\alpha$. We give $X_{\gamma, \gamma^+}$ by the
composite $\sigma_\gamma \circ (t \otimes X_\gamma)$.
\end{itemize}
The following is now Theorem 23.3 of \cite{Ke80}.
\begin{Prop}
Let $\V$ be a cocomplete monoidal category in which each functor $(\thg)
\otimes X \colon \V \to \V$ preserves connected colimits; and let $t \colon I
\to T$ be a pointed object of $\V$. If the free monoid sequence for $(T, t)$
converges at stage $\gamma$, then $X_\gamma$ is the free monoid on $(T, t)$,
with the universal map given by $X_{1, \gamma} \colon T \to X_\gamma$.
\end{Prop}

In fact, this result is a mild generalisation of \cite{Ke80}, since we require
$(\thg) \otimes X$ to preserve only connected colimits, rather than all
colimits; but it is trivial to check that this does not affect the argument in
any way.

In order to apply this result, we observe that:

\begin{Prop}
    If $\C$ is a cocomplete category, then $\cat{LNWFS}(\C)$ is also
    cocomplete; and moreover, each functor $(\thg) \otimes X \colon \cat{LNWFS}(\C) \to
    \cat{LNWFS}(\C)$ preserves connected colimits.
\end{Prop}
\begin{proof}
We first note that the category $\cat{FF}(\C)$ may be obtained by taking the
category $[\C^\mathbf 2, \C]$, slicing this over the object $\cod \colon
\C^\mathbf 2 \to \C$; and then coslicing this under the object $\upsilon \colon
\dom \Rightarrow \cod$ given by $\upsilon_f = f$ for all $f \in \C^\mathbf 2$.
Consequently, $\cat{FF}(\C)$ will be cocomplete whenever $\C$ is. But by
Theorem \ref{diamondtwofold}, the functor $U \colon \cat{LNWFS}(\C) \to
\cat{FF}(\C)$ is a forgetful functor from a category of comonoids, and as such
creates colimits, so that $\cat{LNWFS}(\C)$ is also cocomplete.

In order to see that each functor $(\thg) \otimes X \colon \cat{LNWFS}(\C) \to
    \cat{LNWFS}(\C)$ preserves
connected colimits, we consider the composite
\begin{equation*}
    V := \cat{LNWFS}(\C) \xrightarrow {U} \cat{FF}(\C) \xrightarrow {d_0 \circ (\thg)} [\C^\mathbf 2, \C^\mathbf 2]\text,
\end{equation*}
where we recall that postcomposing with $d_0$ sends a functorial factorisation
$F \colon \C^\mathbf 2 \to \C^\mathbf 3$ to the corresponding endofunctor $R
\colon \C^\mathbf 2 \to \C^\mathbf 2$. It is easy to see that $d_0 \circ
(\thg)$ creates connected colimits; and since $U$ creates all colimits, we
conclude that $V$ creates connected colimits.

Now observe that $V$ sends the monoidal structure on $\cat{LNWFS}(\C)$ to the
compositional monoidal structure on $[\C^\mathbf 2, \C^\mathbf 2]$, so that we
have the following commutative diagram:
\begin{equation*}
\cd{
    \cat{LNWFS}(\C) \ar[r]^{(\thg) \otimes X} \ar[d]_{V} & \cat{LNWFS}(\C) \ar[d]^{V} \\
    [\C^\mathbf 2, \C^\mathbf 2] \ar[r]_{(\thg) \circ VX} & [\C^\mathbf 2, \C^\mathbf 2]\text.
}
\end{equation*}
We wish to show that $(\thg) \otimes X$ preserves connected colimits: but
because $V$ creates them, it suffices to show that the composite around the top
preserves connected colimits; and this follows from the fact that both functors
$V$ and $(\thg) \circ VX$ around the bottom preserve connected colimits.
\end{proof}

Thus the free monoid on $X \in \cat{LNWFS}(\C)$ will exist whenever the free
monoid sequence for $! \colon I \to X$ converges. Sufficient conditions for
convergence are given by Theorem 15.6 of \cite{Ke80}, which when adapted to the
present situation becomes:

\begin{Prop}\label{convergnece}
Let $\V$ be a cocomplete monoidal category, and let $t \colon I \to T$ be a
pointed object of $\V$. If the functor $T \otimes (\thg) \colon \V \to \V$
preserves either $\lambda$-filtered colimits; or $\lambda$-indexed unions of
$\M$-subobjects for some proper, well-copowered $(\E, \M)$ on $\V$, then the
free monoid sequence for $(T, t)$ converges.
\end{Prop}

\pgph There is a problem if we apply this result with $\V = \cat{LNWFS}(\C)$,
since the second of the two smallness criteria requires a proper,
well-copowered $(\E, \M)$ on $\V$; and even if we have such an $(\E, \M)$ on
the category $\C$, we will not, in general, be able to lift it to
$\cat{LNWFS}(\C)$. In order to resolve this problem, we consider again the
composite
\begin{equation*}
    V := \cat{LNWFS}(\C) \xrightarrow {U} \cat{FF}(\C) \xrightarrow {d_0 \circ (\thg)} [\C^\mathbf 2, \C^\mathbf 2]\text.
\end{equation*}

We saw above that this preserves both connected colimits and monoidal
structure; and so takes the free monoid sequence on $! \colon I \to X$ in
$\cat{LNWFS}(\C)$ to the \emph{free monad sequence} on the underlying pointed
endofunctor $\Lambda \colon \id_{\C^\mathbf 2} \Rightarrow R$ of $X$. Moreover,
$V$ reflects isomorphisms: hence the convergence of the latter sequence
guarantees the convergence of the former.

Thus, it will suffice to apply Proposition \ref{convergnece} for $\V =
[\C^\mathbf 2, \C^\mathbf 2]$, which avoids the problem described above, since
any proper, well-copowered $(\E, \M)$ on $\C$ can be lifted without trouble to
$[\C^\mathbf 2, \C^\mathbf 2]$. In fact, it will suffice to lift to $\C^\mathbf
2$, since when we instantiate Proposition \ref{convergnece} at $\V =
[\C^\mathbf 2, \C^\mathbf 2]$, the requirement that $T \otimes (\thg) \colon
[\C^\mathbf 2, \C^\mathbf 2] \to [\C^\mathbf 2, \C^\mathbf 2]$ should preserve
$\lambda$-filtered colimits or unions may be safely reduced to the requirement
that $T \colon \C^\mathbf 2 \to \C^\mathbf 2$ should preserve the same.

We may summarise this argument as follows:

\begin{Prop}\label{reflection2}
Let there be given a cocomplete category $\C$; and let $(F, \Ll) \in
\cat{LNWFS}(\C)$. If the functor $R = d_0 \circ F \colon \C^\mathbf 2 \to
\C^\mathbf 2$ preserves either $\lambda$-filtered colimits; or
$\lambda$-indexed unions of $\M$-subobjects for some proper, well-copowered
$(\E, \M)$ on $\C^\mathbf 2$, then the free monoid sequence for $! \colon I \to
(F, \Ll)$ converges: and in particular, the reflection of $(F, \Ll)$ along
$\G_1 \colon \cat{NWFS}(\C) \to \cat{LNWFS}(\C)$ exists.
\end{Prop}

We are now ready to prove the first part of our main theorem:

\begin{Prop}
Let $\C$ be a cocomplete category satisfying one of the smallness conditions
(*) or (\dag), and let $I \colon \J \to \C^\mathbf 2$ be a category over
$\C^\mathbf 2$ with $\J$ small. Then the free n.w.f.s.\ on $\J$ exists.
\end{Prop}
\begin{proof}
By Proposition \ref{Kanext} and Proposition \ref{reflection}, we may find an
object $(F, \Ll) \in \cat{LNWFS}(\C)$ which is a reflection of $\J$ along $\G_3
\G_2 \colon \cat{LNWFS}(\C) \to \cat{CAT}/\C^\mathbf 2$. We now wish to apply
Proposition \ref{reflection2} to $(F, \Ll)$: so for a $\C$ satisfying (*), we
will show that $R = d_0 \circ F \colon \C^\mathbf 2 \to \C^\mathbf 2$ preserves
$\lambda$-filtered colimits for some $\lambda$; whilst for a $\C$ satisfying
(\dag), we will show that $R$ preserves $\lambda$-indexed unions of
$\M$-subobjects for the induced factorisation system $(\E, \M)$ on $\C^\mathbf
2$. Since the proof is the same in both cases, we restrict our attention to the
former.

We begin by considering the following diagram:
\begin{equation*}
\cd[@R-2em]{ & & \C^\mathbf 2 \\
    \C^\mathbf 2 \ar[r]^{F} & \C^\mathbf 3 \ar[ur]^{d_0} \ar[dr]_{d_2} \\ & & \C^\mathbf 2
}
\end{equation*}
The upper composite is $R \colon \C^\mathbf 2 \to \C^\mathbf 2$, which we are
to show preserves $\lambda$-filtered colimits; but since $d_0$ and $d_2$
preserve and reflect connected colimits, we may equally well show that the
lower composite $L = d_2 \circ F$ preserves $\lambda$-filtered colimits.

Now, from Proposition \ref{Kanext} and Proposition \ref{reflection}, the
functor $L$ has the following explicit description. First we form the left Kan
extension of $I \colon \J \to \C^\mathbf 2$ along itself to obtain a functor $M
\colon \C^\mathbf 2 \to \C^\mathbf 2$. We may describe this by the usual coend
formula
\begin{equation*}
    M(f) = \int^{j \in \J} \C^\mathbf 2(I(j), f) \cdot I(j)\text.
\end{equation*}

We now consider the counit transformation $\epsilon \colon M \Rightarrow
\id_{\C^\mathbf 2}$, whose component at $f$ is the map
\begin{equation*}
    \epsilon_f \colon \int^{j \in \J} \C^\mathbf 2(I(j), f) \cdot I(j) \to f
\end{equation*}
corresponding to the identity transformation $\C^\mathbf 2(I(\thg), f)
\Rightarrow \C^\mathbf 2(I(\thg), f)$; and we factor this transformation
$\epsilon$ as
\begin{equation*}
    M \xRightarrow \xi L \xRightarrow {\Phi} \id_{\C^\mathbf 2}\text,
\end{equation*}
where each component of $\xi$ is a pushout; and each component of $\Phi$ is the
identity in its domain.

Let us first show that $L$ preserves any colimit which $M$ does. Suppose that
$A \colon \I \to \C^\mathbf 2$ is a small diagram whose colimit is preserved by
$K$, and consider the following diagram:
\begin{equation}\label{rendering}
    \cd[@C+2em]{
    \colim_i MA_i \ar[r]^{\colim_i \xi_{A_i}} \ar[d]_{\textsf{can}_M} &
    \colim_i LA_i \ar[r]^{\colim_i \Phi_{A_i}} &
    \colim_i A_i \ar[d]^{=} \\
    M \colim_i A_i \ar[r]_{\xi_{\colim_i A_i}}&
    L \colim_i A_i \ar[r]_{\Phi_{\colim_i A_i}} &
    \colim_i A_i\text.
    }
\end{equation}
The class $\P$ of morphisms in $\C^\mathbf 2$ which are pushout squares is the
left class of a strong factorisation system, and hence stable under colimit:
and thus not only $\xi_{\colim_i A_i}$, but also $\colim_i \xi_{A_i}$, is in
$\P$. Likewise, the class $\D$ of morphisms in $\C^\mathbf 2$ which are
domain-isomorphisms is also the left class of a strong factorisation system on
$\C^\mathbf 2$, whose corresponding right class is the class of
codomain-isomorphisms. Hence $\D$ is also stable under colimit; and so both
$\Phi_{\colim_i A_i}$ and $\colim_i \Phi_{A_i}$ are in $\D$.

The orthogonality property for $(\P, \D)$ now implies that there is a unique
map $\phi \colon \colim_i LA_i \to L \colim A_i$ rendering \eqref{rendering}
commutative; and moreover, that $\phi$ is invertible, since $\textsf{can}_M$
is. But the canonical morphism $\textsf{can}_L \colon \colim_i LA_i \to L
\colim A_i$ makes \eqref{rendering} commute; and so we deduce that
$\textsf{can}_L = \phi$ is invertible as required.

Thus $L$ preserves any colimit which $M$ does: so we will be done if we can
find some $\lambda$ for which $M$ preserves $\lambda$-filtered colimits. Now,
for each $j \in \J$, we have the morphism $I(j) \colon X \to Y$ of $\C^\mathbf
2$: and by condition (*), we can find a $\lambda_j$ for which both $X$ and $Y$
are $\lambda_j$-presentable; from which it follows that $I(j)$ is
$\lambda_j$-presentable in $\C^\mathbf 2$. Thus, if we take $\lambda$ to be a
regular cardinal larger than each $\lambda_j$, then each $I(j)$ is
$\lambda$-presentable in $\C^\mathbf 2$.

We now show that $K$ preserves $\lambda$-filtered colimits. Indeed, suppose
that $A \colon \I \to \C^\mathbf 2$ is a $\lambda$-filtered diagram in
$\C^\mathbf 2$; then we have that
\begin{align*}
    M(\colim_i A_i) & = \textstyle\int^{j \in \J} \C^\mathbf 2(I(j), \colim_i A_i) \cdot I(j) \\
                  & \cong \textstyle\int^{j \in \J} (\colim_i \C^\mathbf 2(I(j), A_i)) \cdot I(j) \text{\ \ \  (as $I(j)$ is $\lambda$-presentable)}\\
                  & \cong \colim_i \textstyle\int^{j \in \J} \C^\mathbf 2(I(j), A_i) \cdot I(j) \text{\ \ \  (as colimits commute with colimits)}\\
                  & = \colim_i M(A_i)\text,
\end{align*}
as desired.
\end{proof}

\section{Constructively-free implies algebraically-free}\label{Sec:cfaf}
In this Section, we prove that all free n.w.f.s.'s obtained by the procedure of
the previous Section are algebraically-free. In order to do this, we will need
to establish a link between our notion of algebraically-free n.w.f.s., and
\cite{Ke80}'s notion of algebraically-free monad. We begin, therefore, by
recalling the latter.

\pgph Let $\sigma \colon \id \Rightarrow S$ be a pointed endofunctor on some
category $\V$. An \emph{$S$-algebra} is an object $X \in \V$ together with a
morphism $x \colon SX \to X$ satisfying $x . \sigma = \id_X$; and an
\emph{$S$-algebra morphism} $(X, x) \to (Y, y)$ is a morphism $f \colon X \to
Y$ of $\V$ for which $f . x = y . Sf$. We write $S$-$\cat{Alg}$ for the
category of $S$-algebras and $S$-algebra morphisms. A \emph{morphism of pointed
endofunctors} $(S, \sigma) \Rightarrow (T, \tau)$ is a natural transformation
$\alpha \colon S \Rightarrow T$ satisfying $\tau = \alpha . \sigma$; and any
such morphism induces a functor $\alpha^\ast \colon T\text-\cat{Alg} \to
S\text-\cat{Alg}$ sending $(X, x)$ to $(X, x . \alpha_X)$.

If we are given a monad $\mathsf T = (T, \eta, \mu)$ on $\V$, we can consider
its category $\mathsf T$-$\cat{Alg}$ of algebras \emph{qua} monad; or we can
consider its category $T$-$\cat{Alg}$ of algebras \emph{qua} pointed
endofunctor. Evidently, every $\mathsf T$-algebra is a $T$-algebra, and so we
have an inclusion functor $\textsf{inc} \colon \mathsf T\text-\cat{Alg} \to
T\text-\cat{Alg}$.

Now let $(S, \sigma)$ be a pointed endofunctor on $\V$. We say that a monad
$\mathsf T = (T, \eta, \mu)$ is \emph{algebraically-free} on $(S, \sigma)$ if
we can provide a morphism of pointed endofunctors $\alpha \colon (S, \sigma)
\Rightarrow (T, \eta)$ such that the composite
\begin{equation*}
    \mathsf T\text-\cat{Alg} \xrightarrow{\textsf{inc}} T\text-\cat{Alg} \xrightarrow{\alpha^\ast} S\text-\cat{Alg}
\end{equation*}
is an isomorphism of categories.

The main result we will need about algebraically-free monads is the following,
which is Theorem 22.3 of \cite{Ke80}:

\begin{Prop}\label{kellyconstr}
Let $\V$ be a cocomplete category, and let $(S, \sigma)$ be a pointed
endofunctor of $\V$. If the free monad sequence $X \colon \cat{On} \to [\V,
\V]$ for $(S, \sigma)$ converges at stage $\gamma$, then the morphism $X_{1,
\gamma} \colon S \Rightarrow X_\gamma$ exhibiting $X_\gamma$ as the free monad
on $S$ also exhibits it as the algebraically-free monad on $S$.
\end{Prop}

\pgph We are now ready to prove the second part of our main Theorem. We suppose
given a cocomplete $\C$, so that any small $I \colon \J \to \C^\mathbf 2$ over
$\C^\mathbf 2$ has a reflection $(F', \Ll')$ along $\G_3 \G_2$; and we now say
that the free n.w.f.s.\ on such a $\J$ \emph{exists constructively} just when
the free monoid sequence for $(F', \Ll')$ converges.

\begin{Prop}\label{constralg}
Let $\C$ be a cocomplete category, and let $I \colon \J \to \C^\mathbf 2$ be a
small category over $\C^\mathbf 2$. If the free n.w.f.s.\ on $\J$ exists
constructively, then it is algebraically-free on $\J$.
\end{Prop}

\begin{proof}
Let us write $(\Ll, \Rr)$ for the free n.w.f.s.\ on $\J$, and $(F', \Ll')$ for
the reflection of $\J$ along $\G_3 \G_2$. By constructive existence, we obtain
$(\Ll, \Rr)$ as the convergent value $X_\gamma$ of the free monoid sequence on
$(F', \Ll')$; and so if $\eta \colon \J \to \Lmap$ exhibits $(\Ll, \Rr)$ as the
free n.w.f.s.\ on $\J$, then the corresponding morphism $\alpha \colon (F',
\Ll') \to (F, \Ll)$ of $\cat{LNWFS}(\C)$ is the map $X_{1, \gamma}$ of this
free monoid sequence. Now, applying the functor
\begin{equation*}
    V:= \cat{LNWFS}(\C) \xrightarrow{U} \cat{FF}(\C) \xrightarrow{d_0 \circ (\thg)} [\C^\mathbf 2, \C^\mathbf 2]
\end{equation*}
to this free monoid sequence yields the free monad sequence for the pointed
endofunctor $\Lambda' \colon \id_{\C^\mathbf 2} \Rightarrow R'$: and the
convergence of the former guarantees the convergence of the latter. Thus by
Proposition \ref{kellyconstr}, we deduce that the map of pointed endofunctors
    $\alpha_r \colon (R', \Lambda') \to (R, \Lambda)\text,$
obtained by applying $V$ to $\alpha$, exhibits $\Rr$ as the algebraically-free
monad on $(R', \Lambda')$.

We now consider the following diagram:
\begin{equation}\label{twosquares}
\cd{
    \mathsf R\text-\cat{Map} \ar[r]^{\textsf{\upshape lift}} \ar[d]_{\id} & \Ll\text-\cat{Map}^{\pitchfork} \ar[r]^-{\eta^\pitchfork} \ar@{.>}[d]^G & \J^\pitchfork \ar@{.>}[d]^H \\
    \mathsf R\text-\cat{Alg} \ar[r]_{\textsf{\upshape inc}} & R\text-\cat{Alg} \ar[r]_{(\alpha_r)^\ast} & R'\text-\cat{Alg}\text.
}
\end{equation}
By algebraic-freeness of $\Rr$, the composite along the bottom is an
isomorphism; and we would like to deduce that the composite along the top is an
isomorphism. To do this, it suffices to find isomorphisms $G$ and $H$ as
indicated which make both squares commute.

We begin by constructing $G$. Recall that an object of
$\Ll\text-\cat{Map}^{\pitchfork}$ is a pair $(g, \phi)$ consisting of a
morphism $g \colon C \to D$ and a mapping $\phi$ which to each object $a \in
\Lmap$ and square
\begin{equation*}
    \cd{
        A
            \ar[r]^-h
            \ar[d]_{U_\Ll(a)} &
        C
            \ar[d]^g \\
        B
            \ar[r]_-k &
        D
    }
\end{equation*} in $\C$,
assigns a fill-in $\phi(a, h, k) \colon B \to C$ which is natural with respect
to morphisms of $\Lmap$. Now, to give such a $\phi$ is equally well to give a
natural transformation
\begin{equation*}
    \phi \colon \C^\mathbf 2(U_\Ll(\thg), g) \Rightarrow \C^\mathbf 2(U_\Ll(\thg), \id_C) \colon (\Lmap)^\op \to \cat{Set}
\end{equation*}
which is a section of the natural transformation
    $\C^\mathbf 2(U_\Ll(\thg), \id_C) \Rightarrow \C^\mathbf 2(U_\Ll(\thg), g)$
induced by postcomposition with $(\id_C, g) \colon \id_C \to g$. But $U_\Ll
\colon \Lmap \to \C^\mathbf 2$ has a right adjoint given by the cofree functor
$C_\Ll \colon \C^\mathbf 2 \to \Lmap$; and thus we have an isomorphism
\begin{equation*}
    \C^\mathbf 2(U_\Ll(\thg), g) \cong \Lmap(\thg, C_\Ll(g))\text.
\end{equation*}
So $\C^\mathbf 2(U_\Ll(\thg), g)$ is represented by $C_\Ll(g)$; and thus by the
Yoneda Lemma, $\phi$ is uniquely determined by where it sends the counit map
$U_\Ll C_\Ll g \to g$; which is to say, by the fill-in it provides for the
square
\begin{equation*}
    \cd{
        C
            \ar[r]^-{\id}
            \ar[d]_{\lambda_g} &
        C
            \ar[d]^g \\
        Kf
            \ar[r]_-{\rho_g} &
        D
    }
\end{equation*} in $\C$. But to provide a fill-in for this square is precisely to make $g$ into an algebra for the pointed endofunctor $(R, \Lambda)$.
Thus we have an isomorphism between objects of $\Lmap^\pitchfork$ and objects
of $R\text-\cat{Alg}$; and it is now straightforward to extend this to the
required isomorphism of categories $G \colon \Lmap^\pitchfork \to
R\text-\cat{Alg}$, and to verify that this $G$ makes the left-hand square of
\eqref{twosquares} commute.

We now complete the proof by constructing the isomorphism $H \colon
\J^\pitchfork \to R'\text-\cat{Alg}$. Proceeding as above, we see that to give
an object of $\J^\pitchfork$ is to give a morphism $g \colon C \to D$ of $\C$
together with a natural transformation
    $\phi \colon \C^\mathbf 2(I(\thg), g) \Rightarrow \C^\mathbf 2(I(\thg), \id_C)$
which is a section of the natural transformation $\C^\mathbf 2(I(\thg), \id_C)
\Rightarrow \C^\mathbf 2(I(\thg), g)$ induced by postcomposing with $(\id_C, g)
\colon \id_C \to g$. Now, if we write
\begin{equation*}
    Mg := \int^{j \in \J} \C^\mathbf 2(I(j), g) \cdot I(j)
\end{equation*}
and $\epsilon_g$ for the counit map $Mg \to g$ as before, then to give $\phi$
is equivalently to give a morphism $k \colon Mg \to \id_C$ satisfying
$\epsilon_g = (\id_C, g) \circ k$. Furthermore, we obtain $L'g$ from $Mg$ by
factorising $\epsilon_g$ as
\begin{equation*}
    \epsilon_g = Mg \xrightarrow{\xi_g} L'g \xrightarrow{\Phi'_g} g\text,
\end{equation*}
where $\xi_g$ is a pushout square, and $\Phi'_g$ is the identity in its domain;
and so given such a map $k \colon Mg \to \id_C$, applying unique
diagonalisation to the diagram
\begin{equation*}
\cd{
    Mg \ar[r]^{k} \ar[d]_{\xi_g} & \id_C \ar[d]^{(\id_C, g)} \\
    L'g \ar[r]_{\Phi'_g} & g
}
\end{equation*}
shows that $k$ is induced by a unique morphism $m \colon L'g \to \id_C$. But to
give such a morphism is to give a diagonal fill-in for the square
\begin{equation*}
    \cd{
        C
            \ar[r]^-{\id}
            \ar[d]_{\lambda'_g} &
        C
            \ar[d]^g \\
        K'f
            \ar[r]_-{\rho'_g} &
        D
    }
\end{equation*} in $\C$; which in turn is to make $g$ into an algebra for the pointed endofunctor $(R', \Lambda')$. The remaining details are again straightforward.
\end{proof}

\section{Comparison with the small object argument}
Since we have advertised the argument of Theorem \ref{mainthm} as an adaptation
of the small object argument, it behooves us to investigate the relationship
between the two. To do this, we combine our main Theorem with Proposition
\ref{relateback} to deduce:
\begin{Prop}\label{compare}
Let $\C$ be a cocomplete category satisfying either of the smallness conditions
(*) or (\dag); and let $J$ be a set of maps in $\C$. Then the w.f.s.\ $(\ELL,
\R)$ cofibrantly generated by $J$ exists.
\end{Prop}

\pgph Since the two classes of maps $\ELL$ and $\R$ of this w.f.s.\ are
entirely determined by the equations $\ELL = {}^{\pitchfork}\R$ and $\R =
J^\pitchfork$, the content of this Proposition is that we may find an $(\ELL,
\R)$-factorisation for every map of $\C$. This is also the content of the small
object argument, and so we may compare the two by comparing the choices of
factorisation which they provide. For a detailed account of the small object
argument, we refer the reader to \cite{Bous} or \cite{Hovey}.

\pgph Suppose we are given a category $\C$ and a set of maps $J$ as in the
Proposition; and let $g \colon C \to D$ be a morphism of $\C$ that we wish to
factorise. The first step in both the small object argument and our argument
turns out to be the same. In the small object argument, we form the set $S$
whose elements are squares
\begin{equation*}
\cd{
    A \ar[r]^{h} \ar[d]_f & C \ar[d]^g \\
    B \ar[r]_k & D
}
\end{equation*}
such that $f \in J$. We then form the coproduct
\begin{equation}\label{obvioussquare}
\cd[@C+1.5em]{
    \sum_{x \in S} A_x \ar[r]^-{[h_x]_{x \in S}} \ar[d]_{\sum_{x \in S} f_x} & C \ar[d]^g \\
    \sum_{x \in S} B_x \ar[r]_-{[k_x]_{x \in S}} & D
}
\end{equation}
and define an object $K'g$ and morphisms $\lambda'_g$ and $\rho'_g$ by
factorising this square as
\begin{equation}\label{pushoutsquare}
\cd[@C+1.5em]{
    \sum_{x \in S} A_x \ar[r]^-{[h_x]_{x \in S}} \ar[d]_{\sum_{x \in S} f_x} & C \ar[d]^{\lambda'_g} \ar[r]^{\id_C} & C \ar[d]^g \\
    \sum_{x \in S} B_x \ar[r]_-{\xi_g} & K'g \ar[r]_{\rho'_g} & D\text,
}
\end{equation}
where the left-hand square is a pushout.

On the other hand, suppose we view $J$ as a discrete subcategory $\J$ of
$\C^\mathbf 2$; and write $I \colon \J \hookrightarrow \C^\mathbf 2$ for the
inclusion functor. Then we may view \eqref{obvioussquare} as the morphism
\begin{equation*}
    \epsilon_g \colon \int^{f \in \J} \C^\mathbf 2(If, g) \cdot If \to g\text,
\end{equation*}
of $\C^\mathbf 2$; which is to say, the component at $g$ of the counit
transformation
\begin{equation*}
    \epsilon \colon \Lan_I(I) \Rightarrow \id_{\C^\mathbf 2} \colon \C^\mathbf 2 \to \C^\mathbf 2\text.
\end{equation*}
We may then view \eqref{pushoutsquare} as the component at $g$ of the
factorisation of $\epsilon$ into a map which is componentwise a pushout,
followed by a map whose domain components are identities. Thus the assignation
$g \mapsto (\lambda'_g, \rho'_g)$ obtained from the small object argument is
just the underlying factorisation of the reflection of $I \colon \J
\hookrightarrow \C^\mathbf 2$ along $\G_3 \G_2$.

\pgph At this point, the two arguments under consideration diverge from each
other. The small object argument is the more naive of the two: it simply
iterates the above procedure, each time replacing the map $g$ with the map
$\rho'_g$. This gives rise to the countable sequence
\begin{equation*}
\cd{
    C \ar[d]_g \ar[r]^{\lambda'_g} & K'g \ar[d]_{\rho'_g} \ar[r]^{\lambda'_{\rho'_g}} & K'\rho'_g \ar[d]_{\rho'_{\rho'_g}} \ar[r]^{\lambda'_{\rho'_{\rho'_g}}} & \dots \\
    D \ar[r]_{\id_D} & D \ar[r]_{\id_D} & D \ar[r]_{\id_D} & \dots\text,
}
\end{equation*}
which we extend transfinitely by taking colimits at limit ordinals. However, as
pointed out in \cite{injective}, this sequence \emph{almost never converges}.
Instead, the small object argument requires one to choose an arbitrary ordinal
at which to stop: or rather, an ordinal which is large enough to ensure that
the right part of the corresponding factorisation lies in $J^\pitchfork$.

\pgph Our argument produces a different transfinite sequence, whose first few
terms are:
\begin{equation*}
\cd{
    C \ar[d]_g \ar[r]^{\lambda'_g} & K'g \ar[d]_{\rho'_g} \ar[r]^{\lambda''_g} & K''g \ar[d]_{\rho''_g} \ar[r]^{\lambda'''_g} & \dots \\
    D \ar[r]_{\id_D} & D \ar[r]_{\id_D} & D \ar[r]_{\id_D} & \dots\text;
}
\end{equation*}
here, $K''g$ is the coequaliser
\begin{equation*}
\cd[@C+3em]{
    K'g \ar@<4pt>[r]^{\lambda'_{\rho'_g}} \ar@<-4pt>[r]_{K'(\lambda'_g, \id_D)} & K'\rho'_g \ar[r] & K''g\text,
}
\end{equation*}
and in general, the term at stage $\alpha$ in this sequence will be a quotient
of the corresponding term at stage $\alpha$ in the small object argument.

We may understand this quotienting process as follows. The small object
argument provides a way of taking a map $g \colon C \to D$, and recursively
adding elements to its domain which witness the required lifting properties
against the set $J$. This process must be recursive, since the process of
adding witnesses can create new instances of the lifting properties: which in
turn will require new witnesses to be added, and so on.

 However, the small
object argument is badly behaved: at each stage it adds new witnesses for
\emph{all} instances of the required lifting properties~--~including those
instances for which witnesses were added at a previous stage of the induction.
The effect of the quotienting process which our argument carries out is to
collapse these superfluous new witnesses back onto their predecessors.

\section{Applications}
We end the paper with two simple applications of Theorem \ref{mainthm}.

\pgph In Examples \ref{exs1}, we saw that the set $J$ of horn inclusions
generates a plain w.f.s.\ (anodyne extensions, Kan fibrations) on $\cat{SSet}$.
If we view the set $J$ as a discrete subcategory $\J \hookrightarrow
\cat{SSet}^\mathbf 2$, then it also generates a natural w.f.s.\ $(\Ll, \Rr)$.

By restricting the monad $\Rr \colon \cat{SSet}^\mathbf 2 \to
\cat{SSet}^\mathbf 2$ of this natural w.f.s.\ to the slice over the terminal
object, we obtain a monad $T \colon \cat{SSet} \to \cat{SSet}$, whose category
of algebras is the category $\cat{AlgKan}$ of ``algebraic Kan complexes'':
simplicial sets equipped with a chosen filler for every horn, subject to no
further coherence conditions. Since $\cat{AlgKan}$ is finitarily monadic over
$\cat{SSet}$, it is locally finitely presentable, and hence provides a rich
categorical base for further constructions.

Observe that the morphisms of $\cat{AlgKan}$ are maps of simplicial sets which
\emph{strictly} preserve the chosen fillers. Whilst these maps are of some
theoretical importance, we are more likely to be interested in the category
$\cat{AlgKan}_\psi$ whose objects are the same, but whose morphisms are
arbitrary maps of simplicial sets. We may obtain this category by considering
the adjunction
\begin{equation*}
  \cd{\cat{AlgKan} \ar@<4pt>[r]^-{U} \ar@{}[r]|-{\top} & \cat{SSet} \ar@<4pt>[l]^-{F}\text.}
\end{equation*}
This generates a comonad $FU$ on $\cat{AlgKan}$; and the corresponding
co-Kleisli category is precisely $\cat{AlgKan}_\psi$. In particular, we deduce
that the inclusion functor $\cat{AlgKan} \hookrightarrow \cat{AlgKan}_\psi$ has
a left adjoint. It is a corresponding result which forms the cornerstone of
two-dimensional monad theory \cite[Theorem 3.13]{2dmonad}.

\pgph For an example even more in the spirit of \cite{2dmonad}, we consider the
category $\C = \cat{2}\text-\cat{Cat}$ and the set of maps $J$ given as
follows:
  \[\cd[@R+2em]{\emptyset \ar@{.>}[d] \\ \bullet}\text; \qquad \cd[@R+2em]{ \bullet \ar@{}[r]_{}="a" & \bullet  \\
  \bullet \ar[r]^{}="b" \ar@{.>}"a"; "b"& \bullet}\text; \qquad \cd[@R+2em]{ \bullet \ar@/^1em/[r] \ar@/_1em/[r]_{}="a" & \bullet  \\
  \bullet \ar@/^1em/[r]^{}="b" \ar@/_1em/[r] \ar@{}[r] \ar@{=>}?(0.5)+/u  0.15cm/;?(0.5)+/d 0.15cm/ \ar@{.>}"a"; "b"& \bullet}\text; \qquad
  \cd[@R+2em]{ \bullet \ar@/^1em/[r] \ar@/_1em/[r]_{}="a" \ar@{}[r] \ar@{=>}?(0.65)+/u  0.15cm/;?(0.65)+/d 0.15cm/ \ar@{}[r] \ar@{=>}?(0.35)+/u  0.15cm/;?(0.35)+/d 0.15cm/ & \bullet  \\
  \bullet \ar@/^1em/[r]^{}="b" \ar@/_1em/[r] \ar@{}[r] \ar@{=>}?(0.5)+/u  0.15cm/;?(0.5)+/d 0.15cm/ \ar@{.>}"a"; "b"& \bullet\text.}\]

These maps generate a plain w.f.s.\ which is one half of the model structure on
$\cat{2}\text-\cat{Cat}$ described by Lack in \cite{2catmodel}. Our purpose
here will be to consider the corresponding natural w.f.s.\ $(\Ll, \Rr)$
generated by these maps, where as usual we view $J$ as a discrete subcategory
$\J \hookrightarrow \C^\mathbf 2$.

In particular, if we take the comonad $\Ll$ for this natural w.f.s.\ and
restrict it to the coslice under the initial object, we obtain a comonad $Q
\colon \cat{2}\text-\cat{Cat} \to \cat{2}\text-\cat{Cat}$. We can describe $Q$
quite explicitly. Given a 2-category $\K$, we first form the free 2-category
$FU\K$ on the underlying 1-graph of $\K$. Then we take the counit 2-functor
$\epsilon_\K \colon FU\K \to \K$ and factorise it as
\begin{equation*}
    \epsilon_\K = FU\K \xrightarrow{\xi_\K} Q\K \xrightarrow{\phi_\K} \K
\end{equation*}
where $\xi_\K$ is bijective on objects and 1-cells, and $\phi_\K$ is locally
fully faithful. The resultant $Q\K$ is precisely the ``homomorphism
classifier'' of $\K$: it is characterised by an isomorphism, natural in $\ELL$,
between
\begin{equation*}
    \text{2-functors } Q\K \to \ELL \qquad \text{and} \qquad \text{pseudofunctors } \K \to \ELL\text.
\end{equation*}
It follows from this characterisation that the co-Kleisli category of $Q$ is
the category $\cat{2}\text-\cat{Cat}_\psi$ of 2-categories and pseudofunctors
between them.

Observe that in this example, we at no point had to define what a
``pseudofunctor'' was: it emerged simply from applying our apparatus for a
well-chosen set of maps $J$. Of course, since we already knew what
pseudofunctors were, we did not gain much from this; however, it suggests that
for a more complex $\C$, we may be able to define a suitable notion of
``pseudomorphism'' simply by applying the above argument for a suitable set of
maps $J$.

As an example of this, let us consider the category $\cat{Tricat}$ of
tricategories and (strict) structure-preserving maps between them, and see how
this argument allows us to derive the notion of trihomomorphism. By
``tricategory'', we will mean \cite{nicktricats}'s algebraic definition of
tricategory, so that $\cat{Tricat}$ is finitarily monadic over the category
$\cat{GSet}_3$ of 3-dimensional globular sets; and in particular is locally
finitely presentable. Let us write
\begin{equation*}
  \cd{\cat{Tricat} \ar@<4pt>[r]^-{U} \ar@{}[r]|-{\top} & \cat{GSet}_3 \ar@<4pt>[l]^-{F}}
\end{equation*}
for the free/forgetful adjunction. We define a set $J$ of morphisms in
$\cat{Tricat}$ by taking the following set of maps in $\cat{GSet}_3$:
  \[\cd[@R+2em@C+1em]{\emptyset \ar@{.>}[d] \\ \bullet}\text; \qquad \cd[@R+2em@C+1em]{ \bullet \ar@{}[r]_{}="a" & \bullet  \\
  \bullet \ar[r]^{}="b" \ar@{.>}"a"; "b"& \bullet}\text; \qquad \cd[@R+2em@C+1em]{ \bullet \ar@/^1em/[r] \ar@/_1em/[r]_{}="a" & \bullet  \\
  \bullet \ar@/^1em/[r]^{}="b" \ar@/_1em/[r] \ar@{}[r] \ar@{=>}?(0.5)+/u  0.15cm/;?(0.5)+/d 0.15cm/ \ar@{.>}"a"; "b"& \bullet}\text; \qquad
  \cd[@R+2em@C+1em]{ \bullet \ar@/^1em/[r] \ar@/_1em/[r]_{}="a" \ar@{}[r] \ar@{=>}?(0.65)+/u  0.15cm/;?(0.65)+/d 0.15cm/ \ar@{}[r] \ar@{=>}?(0.35)+/u  0.15cm/;?(0.35)+/d 0.15cm/ & \bullet  \\
  \bullet \ar@{}[r] \ar@3?(0.5)+/l 0.13cm/;?(0.5)+/r 0.13cm/ \ar@/^1em/[r]^{}="b" \ar@/_1em/[r] \ar@{}[r] \ar@{=>}?(0.7)+/u  0.15cm/;?(0.7)+/d 0.15cm/ \ar@{}[r] \ar@{=>}?(0.3)+/u  0.15cm/;?(0.3)+/d 0.15cm/ & \bullet \ar@{.>}"a"; "b"}\text; \qquad
  \cd[@R+2em@C+1em]{ \bullet \ar@{}[r] \ar@3?(0.5)+/l 0.13cm/+/u 0.13cm/;?(0.5)+/r 0.13cm/+/u 0.13cm/ \ar@3?(0.5)+/l 0.13cm/+/d 0.26cm/;?(0.5)+/r 0.13cm/+/d 0.26cm/ \ar@/^1em/[r] \ar@/_1em/[r]_{}="a" \ar@{}[r] \ar@{=>}?(0.7)+/u  0.15cm/;?(0.7)+/d 0.15cm/ \ar@{}[r] \ar@{=>}?(0.3)+/u  0.15cm/;?(0.3)+/d 0.15cm/ & \bullet\\
  \bullet \ar@{}[r] \ar@3?(0.5)+/l 0.13cm/;?(0.5)+/r 0.13cm/ \ar@/^1em/[r]^{}="b" \ar@/_1em/[r] \ar@{}[r] \ar@{=>}?(0.7)+/u  0.15cm/;?(0.7)+/d 0.15cm/ \ar@{}[r] \ar@{=>}?(0.3)+/u  0.15cm/;?(0.3)+/d 0.15cm/ & \bullet\text, \ar@{.>}"a"; "b"}\]
and applying the free functor $F$ to each of them. We now proceed as before: we
consider this set $J$ as a discrete subcategory $\J \hookrightarrow
\cat{Tricat}^\mathbf 2$ and let $(\Ll, \Rr)$ be the n.w.f.s.\ generated by
$\J$; and then let $Q$ be the comonad on $\cat{Tricat}$ given by the
restriction of $\Ll$ to the coslice under the initial object.

We now \emph{define} a trihomomorphism $\mathcal S \to \T$ to be a strict
morphism $Q\mathcal S \to \T$, and define the category $\cat{Tricat}_\psi$ of
tricategories and trihomomorphisms to be the co-Kleisli category of $Q$. The
notion of trihomomorphism we obtain in this way cannot be the one we are used
to from \cite{GPS}, since the latter does not admit a strictly associative
composition: see \cite{gg}. Nonetheless, we can show that our new notion of
trihomomorphism is equivalent to the old one, in that we can exhibit a
biequivalence between a suitably defined 2-category of these new
trihomomorphisms and a corresponding bicategory of the usual ones.

The full details of this will be worked out in a forthcoming paper; but for
now, let us merely say that this method should immediately extend to
(sufficiently algebraic) weak $n$-categories and even weak $\omega$-categories,
thereby allowing us to give a notion of ``weak morphism of
$\omega$-categories'' which admits a strictly associative composition.

\appendix
\section{Algebraically-free implies free}
The purpose of this Appendix is to sketch a proof of the following result:
\begin{Thm}\label{freealgfree}
Let $(\Ll, \Rr)$ be a n.w.f.s.\ on $\C$ which is algebraically-free on $I
\colon \J \to \C^\mathbf 2$. Then $(\Ll, \Rr)$ is free on $\J$.
\end{Thm}
\begin{proof}[Proof]
We first define a monoidal structure on the category $\cat{CAT} / \C^\mathbf
2$. Given $U \colon \A \to \C^\mathbf 2$ and $V \colon \B \to \C^\mathbf 2$,
their tensor product $W \colon \A \otimes \B \to \C^\mathbf 2$ is obtained by
first taking the pullback
\begin{equation*}
    \cd[@-1em]{
 \A \otimes \B
    \ar[rr]
    \ar[dd] \pushoutcorner & &
 \A
    \ar[d]^{U} \\ & &
 \C^\mathbf 2 \ar[d]^{\dom}\\
 \B
    \ar[r]_{V} &
 \C^\mathbf 2
    \ar[r]_{\cod} &
 \C\text,
}
\end{equation*}
and then defining the projection $W$ by $W(a, b) = Ua \circ Vb$. The unit for
this tensor product is the object $(s_0 \colon \C \to \C^\mathbf 2)$, where
$s_0$ is the functor induced by homming the unique map $\sigma_0 \colon \mathbf
2 \to \mathbf 1$ into $\C$; thus $s_0(c) = \id_c \colon c \to c$.

Next, we show that, for any n.w.f.s.\ $(\Ll, \Rr)$ on $\C$, the object $(U_\Ll
\colon \Lmap \to \C^\mathbf 2)$ is a monoid with respect to this monoidal
structure: the key point being that, given $\Ll$-map structures on $f \colon X
\to Y$ and $g \colon Y \to Z$, we may define an $\Ll$-map structure on $gf
\colon X \to Z$. Indeed, if these two $\Ll$-map structures are provided by
morphisms $s \colon Y \to Kf$ and $t \colon Z \to Kg$ (as in \S
\ref{worthdefn}), then the $\Ll$-map structure on the composite $gf$ is given
by:
\begin{equation*}
    Z \xrightarrow{t} Kg \xrightarrow{K(s, \id_Z)} K(g \circ \rho_f) \xrightarrow{K(K(1, g), 1)} K\rho_{gf}
    \xrightarrow{\pi_{gf}} K(gf)\text.
\end{equation*}
The remaining details are routine; and by dualising, we see that $U_\Rr \colon
\Rmap \to \C^\mathbf 2$ is also a monoid in $\cat{CAT}/\C^\mathbf 2$.

We may now show that, if $\alpha \colon (\Ll, \Rr) \to (\Ll', \Rr')$ is a map
of n.w.f.s.'s, then the induced functors $(\alpha_l)_\ast \colon \Lmap \to
\Ll'\text-\cat{Map}$ and $(\alpha_r)^\ast \colon \Rr'\text-\cat{Map} \to \Rmap$
are maps of monoids; so that the semantics functors $\G$ and $\H$ may be lifted
to functors
\begin{align*}
    \hat \G \colon \cat{NWFS}(\C) & \to \cat{Mon}(\cat{CAT}/\C^\mathbf 2) \\
    \text{and} \qquad \hat \H \colon \cat{NWFS}(\C) & \to \big(\cat{Mon}(\cat{CAT}/\C^\mathbf 2)\big)^\op\text.
\end{align*}

We now arrive at a crucial juncture in the proof: we show that $\hat \G$ and
$\hat \H$ are fully faithful. In the case of $\hat \G$, for example, we
consider n.w.f.s.'s $(\Ll, \Rr)$ and $(\Ll', \Rr')$ on $\C$, and a map of
monoids $F \colon \Lmap \to \Ll'\text-\cat{Map}$ over $\C^\mathbf 2$; and must
show that there is a unique morphism $\alpha \colon (\Ll, \Rr) \to (\Ll',
\Rr')$ for which $F = (\alpha_l)_\ast$. To do this, we consider squares of the
following form:
\begin{equation*}
\cd{
    A \ar[d]_{\lambda_f} \ar[r]^{\lambda'_f} &
    K'f \ar[d]^{\rho'_f} \\
    Kf \ar[r]_{\rho_f} &
    B\text.
}
\end{equation*}
We can make $\rho'_f$ into an $\Rr'$-map, since it is the free $\Rr'$-map on
$f$. Similarly, we can make $\lambda_f$ into an $\Ll$-map; and by applying the
functor $F \colon \Lmap \to \Ll'\text-\cat{Map}$, we may make it into an
$\Ll'$-map. Now we apply the lifting operation associated with $(\Ll', \Rr')$
to obtain a morphism $\alpha_f \colon Kf \to K'f$. These maps $\alpha_f$
provide the components of a morphism between the underlying functorial
factorisations of $(\Ll, \Rr)$ and $(\Ll', \Rr')$: it remains only to check
that the comonad and monad structures are preserved. This is just a matter of
checking details, but makes essential use of two facts: that $F$ is a map of
monoids; and that the distributivity axiom holds in $(\Ll, \Rr)$ and $(\Ll',
\Rr')$.

Next, we prove that for any category $U \colon \A \to \C^\mathbf 2$ over
$\C^\mathbf 2$, the category $\A^\pitchfork \to \C^\mathbf 2$ is a monoid in
$\cat{CAT}/\C^\mathbf 2$. The key point is to show that, whenever we equip
morphisms $f \colon C \to D$ and $g \colon D \to E$ of $\C$ with coherent
choices of liftings against the elements of $\A$, we induce a corresponding
equipment on the composite $gf$. Indeed, given $a \in \A$ and a square
\begin{equation*}
    \cd{
        A
            \ar[r]^-h
            \ar[dd]_{Ua} &
        C
            \ar[d]^f \\ &
        D
            \ar[d]^g \\
        B
            \ar[r]_-k &
        E\text,
    }
\end{equation*}
we define $\phi_{gf}(a, h, k) \colon B \to C$ as follows. First we form $j :=
\phi_g(a, fh, k) \colon B \to D$; and now we take $\phi_{gf}(a, h, k) :=
\phi_f(a, h, j) \colon B \to C$.

We may now check that if $F \colon \A \to \B$ is a morphism of
$\cat{CAT}/\C^\mathbf 2$, then the morphism $F^\pitchfork \colon \B^\pitchfork
\to \A^\pitchfork$ respects the monoid structures on $\A^\pitchfork$ and
$\B^\pitchfork$, so that the functors $(\thg)^\pitchfork$, and dually
${}^\pitchfork(\thg)$, lift to functors
\begin{align*}
    (\thg)^\pitchfork & \colon (\cat{CAT}/\C^\mathbf 2)^\op \to \cat{Mon}(\cat{CAT}/\C^\mathbf 2) \\
    \text{and} \qquad {}^\pitchfork(\thg) & \colon \cat{CAT}/\C^\mathbf 2 \to \big(\cat{Mon}(\cat{CAT}/\C^\mathbf 2)\big)^\op\text.
\end{align*}

Finally, we may show that for any n.w.f.s.\ $(\Ll, \Rr)$ on $\C$, the canonical
operation of lifting $\textsf{lift} \colon \Rmap \to \Lmap^\pitchfork$ is a
monoid morphism in $\cat{CAT}/\C^\mathbf 2$. Again, this is simply a matter of
checking details.

We now have all the material we need to prove the Theorem. We suppose ourselves
given a n.w.f.s.\ $(\Ll, \Rr)$ which is algebraically-free on $I \colon \J \to
\C^\mathbf 2$ via the morphism $\eta \colon \J \to \Lmap$: and are required to
show that $(\Ll, \Rr)$ is free on $\J$. So consider a further n.w.f.s.\ $(\Ll',
\Rr')$ on $\C$, and a morphism $F \colon \J \to \Ll'\text-\cat{Map}$ over
$\C^\mathbf 2$. We can form the following diagram of functors over $\C^\mathbf
2$:
\begin{equation*}
    \cd{
      \Rmap \ar[r]^{\textsf{lift}} & \Lmap^\pitchfork \ar[r]^-{\eta^\pitchfork} & \J^\pitchfork\text.\\
      \Rr'\text-\cat{Map} \ar[r]_{\textsf{lift}} & \Ll'\text-\cat{Map}^\pitchfork \ar[ur]_{F^\pitchfork}
    }
\end{equation*}

By algebraic-freeness, the composite along the top is invertible, and so we
obtain from this diagram a functor $\Rr'\text-\cat{Map} \to
\Rr\text-\cat{Map}$. But every map in the diagram is a map of monoids, and
hence the induced functor $\Rr'\text-\cat{Map} \to \Rr\text-\cat{Map}$ is too;
and so is induced by a unique morphism of n.w.f.s.'s $\alpha \colon (\Ll, \Rr)
\to (\Ll', \Rr')$.

It requires a little more work to show $(\alpha_l)_\ast \circ \eta = F$, and
that $\alpha$ is the unique morphism of n.w.f.s.'s with this property. The two
essential facts that we need are that, for any n.w.f.s.\ $(\Ll, \Rr)$, the
canonical morphism $\Lmap \to {}^\pitchfork(\Rmap)$ is a monomorphism; and
that, for any morphism of n.w.f.s.'s $\alpha \colon (\Ll, \Rr) \to (\Ll',
\Rr')$, the following diagram commutes:
\begin{equation*}
    \cd{
    \Rmap \ar[r]^{\textsf{lift}} \ar[d]_{(\alpha_r)^\ast} & \Lmap^\pitchfork \ar[d]^{((\alpha_l)_\ast)^\pitchfork} \\
    \Rr'\text-\cat{Map} \ar[r]_{\textsf{lift}} & \Ll'\text-\cat{Map}^\pitchfork\text.
    }
\end{equation*}
We leave these details to the reader.
\end{proof}

\bibliography{biblio}

\end{document}